\documentclass[12pt,twoside,dvips,draft]{amsart}
\usepackage{amssymb}

\usepackage[all]{xy}
\nonstopmode

\headsep=0.5cm  \numberwithin{equation}{section}
\font\tengothic=eufm10 scaled\magstep 1 \font\sevengothic=eufm7
scaled\magstep 1
\newfam\gothicfam
      \textfont\gothicfam=\tengothic
      \scriptfont\gothicfam=\sevengothic

\newtheorem{theorem}{Theorem}[section]

\newtheorem{proposition}[theorem]{Proposition}
\newtheorem{corollary}[theorem]{Corollary}

\theoremstyle{definition}
\newtheorem{definition}[theorem]{Definition}
\newtheorem{remark}[theorem]{Remark}
\newtheorem{problem}[theorem]{Problem}
\newtheorem{example}[theorem]{Example}

\newtheorem{notation}[theorem]{Notation}
\newtheorem{question}[theorem]{Question}

\newcommand{\Proj}{\operatorname{Proj}}

\newcommand{\cO}{{\mathcal O}}

\newcommand {\RR}{\mathbb{R}}
\newcommand {\ZZ}{\mathbb{Z}}

\newcommand {\PP}{\mathbb{P}}

\newlength\meiabase
\newlength\altura

\newlength\meiabasegrande
\newlength\alturagrande

\newlength\meiabasex
\newlength\alturax
\newlength\meiaalturax

\newlength\meiabasexgrande
\newlength\alturaxgrande

\newlength\meiabasexpeq
\newlength\alturaxpeq

\begin{document}
\title[ Stability and Unobstructedness of Syzygy Bundles]
 {Stability and Unobstructedness of Syzygy Bundles}

\author[L.\ Costa, P. Macias Marques, R.M.\ Mir\'o-Roig]{L.\ Costa$^*$, P. Macias Marques$^{**}$,  R.M.\
Mir\'o-Roig$^{***}$}

\address{Facultat de Matem\`atiques, Departament d'\`Algebra i Geometria, Gran Via de les Corts Catalanes, 585, 08007 Barcelona, SPAIN }

\email{costa@ub.edu, pmm@uevora.pt, miro@ub.edu}

\date{\today}

\thanks{$^*$ Partially supported by MTM2007-61104.}
\thanks{$^{**}$ Partially supported by Funda\c{c}\~ao para a Ci\^encia e Tecnologia, under grant SFRH/BD/27929/2006, and by CIMA -- Centro de Investiga\c{c}\~ao  em Matem\'atica e Aplica\c{c}\~oes, Universidade de \'Evora.}
\thanks{$^{***}$ Partially supported by MTM2007-61104.}



\begin{abstract}
It is a longstanding problem in Algebraic Geometry to determine whether the syzygy bundle $E_{d_1,\ldots ,d_n}$ on $\PP^N$ defined as the kernel of a general epimorphism $\xymatrix{\phi:\cO(-d_1)\oplus\cdots\oplus\cO(-d_n)\ar[r] &\cO}$ is (semi)stable. In this note, we restrict our attention to the case of syzygy bundles $E_{d,n}$ on $\PP^N$ associated to $n$ generic forms $f_1,\ldots,f_n\in K[X_0,X_1,\ldots ,X_N]$ of the same degree $d$. Our first goal is to prove that $E_{d,n}$ is stable if $N+1\le n\le\tbinom{d+2}{2}+N-2$ and ${(N,n,d)\neq(2,5,2)}$. This bound improves, in general, the bound $n\le d(N+1)$ given by G. Hein in \cite{B}, Appendix~A.

In the last part of the paper, we study moduli spaces of stable rank $n-1$ vector bundles on $\PP^N$ containing syzygy bundles. We prove that if $N+1\le n\le\tbinom{d+2}{2}+N-2$, $N\ne 3$ and ${(N,n,d)\neq(2,5,2)}$, then the syzygy bundle  $E_{d,n}$ is  unobstructed and it belongs to a generically smooth irreducible component of dimension $n\tbinom{d+N}{N}-n^2$, if $N \geq 4$, and $n\tbinom{d+2}{2}+n\tbinom{d-1}{2}-n^2$, if $N=2$.
\end{abstract}


\maketitle

\tableofcontents


\section{Introduction}\label{intro}

Let $R=K[X_0,X_1,\ldots,X_N]$, $\PP^N=\Proj(R)$ be the $N$-dimensional projective space over an algebraically closed field $K$ of characteristic 0. Set $\frak{m}=(X_0,X_1,\ldots,X_N)$. It is a classical and difficult problem in Algebraic Geometry, as well in Commutative Algebra, to understand the syzygy bundle $E_{d_1,\ldots ,d_n}$ on $\PP^N$ defined as the kernel of a general epimorphism
\[
\xymatrix{
\phi=(f_1,\ldots ,f_n):
    \cO_{\PP^N}(-d_1)\oplus\cdots\oplus\cO_{\PP^N}(-d_n)
    \ar[r]
    &\cO_{\PP^N},
}
\]
where $(f_1,\ldots,f_n)\subset R$ is an $\frak{m}$-primary ideal, and $f_{i}$ is an homogeneous polynomial of degree $d_i=deg(f_i)$. We would like to know the cohomology of  $E_{d_1,\ldots,d_n}$, its splitting type on a generic line, and whether it is simple, exceptional or stable. In particular, we are led to consider the following problem:

\begin{problem}\label{mainpblm} Let $f_1,\ldots,f_n\in R$ be a family of $\frak{m}$-primary homogeneous polynomials of degree $deg(f_i)=d_i$, $1\le i \le n$. Let $E_{d_1,\ldots ,d_n}$ be the syzygy bundle on $\PP^N$ associated to $f_1,\ldots ,f_n$. Is  $E_{d_1,\ldots ,d_n}$ a (semi)stable vector bundle on $\PP^N$?
\end{problem}

In the last few years, Problem \ref{mainpblm} has been extensively studied  and surprisingly only a few partial results have been obtained. We refer to \cite{B} and \cite{B1} for precise information. In this paper, we restrict our attention to the case $d_1=d_2=\ldots = d_n=d$ and we address the following problem, which should be viewed as a particular case of Problem \ref{mainpblm}.

\begin{problem}\label{pblm} Let $f_1,\ldots,f_n\in R$ be a family of \mbox{$\frak{m}$-primary} forms of
the same degree $d$ and let $E_{d,n}$ be the syzygy bundle
associated to them. Is $E_{d,n}$ a (semi)stable vector bundle on
$\PP^N$?
\end{problem}

Note that since $(f_1,\ldots,f_n)$ is an \mbox{$\frak{m}$-primary} ideal, we always have $N+1\le n\le\tbinom{d+N}{N}$. Problem~\ref{pblm} turns out to be true for a set of $n$ general \mbox{$\frak{m}$-primary} forms of the same degree $d$, provided
\begin{itemize}
\item $d$ and $N$ are arbitrary and $n=\tbinom{N+d}{N}$ \cite{P};
\item $d$ and $N$ are arbitrary  and $n=N+1$ \cite{BS};
\item  $d$ and $N$ are arbitrary and $n\le d(N+1)$ \cite{B}.
\end{itemize}

The first goal of this paper is to give an affirmative answer to Problem~\ref{pblm} for the case of $n$ general \mbox{$\frak{m}$-primary} forms of the same degree $d$, provided
\begin{itemize}
\item[(1)] $N=2$ and $3\le n\le\tbinom{d+2}{2}$ (see Theorem~\ref{mainthm});
\item [(2)] $N\ge 2$ and $N+1\le n\le\tbinom{d+2}{2}+N-2$ (see Theorem \ref{mainthm2}).
\end{itemize}

\vskip 2mm
We want to point out that the result (1) was announced by Brenner in \cite{B} but no proof was included and the result (2) strongly improves, in general,  the bound $N+1\le n \le d(N+1)$ given by G. Hein in \cite{B1}, Theorem A1.

\vskip 2mm

In the last section of this work, we also study the
unobstructedness of stable syzygy bundles on $\PP^N$. There
exists a beautiful theorem due to Maruyama establishing the
existence of the moduli space $M=M(r;c_1, \ldots, c_{s})$  of rank $r$, stable vector bundles $E$ on $\PP^N$ with fixed Chern classes $c_i(E)=c_i$ for $i=1, \ldots, s=min(r,N)$ (see \cite{MA} and \cite{MA1}). Unfortunately, in general, very little is known about its local and global structure. In this paper, we prove that  points $[E_{d,n}]$ of $M=M(r;c_1, \ldots, c_{s})$ parameterizing stable syzygy bundles $E_{d,n}$ on $\PP^N$, $N \neq 3$ and $N+1\le n\le\tbinom{d+2}{2}+N-2$, are smooth and we compute the dimension of the irreducible component of $M=M(r;c_1, \ldots, c_{s})$ passing through $[E_{d,n}]$ in terms of $d$, $n$ and $N$ (see Theorem \ref{moduli}).

\vskip 2mm  \noindent {\bf Notation:} We work over an
algebraically closed field $K$ of characteristic zero. We set
$\PP^N=\Proj(K[X_0,X_1,\ldots ,X_N])$ and $\mathfrak{m}=(X_0,X_1,\ldots ,X_N)$. Given coherent sheaves $E$ and $F$ on $\PP^N$, we write $h^i(E)$ (resp.\ $ext^i(E,F)$) to denote the dimension of the $i$th cohomology group $H^i(\PP^N,E)=H^i(E)$ (resp.\ $i$th Ext group $Ext^i(E,F)$) as a \mbox{$K$-vector} space.

For any $x\in \RR$, we set $\lceil x\rceil:=\min\{n\in \ZZ \mid x\le n \}$.

\section{Stability of syzygy bundles. Generalities}

In this section, we recall the notion of (semi)stability of torsion free sheaves on projective spaces and its basic properties. We review the useful cohomological characterization of (semi)stability due to Hoppe as well as  its applications to the problem of determining the (semi)stability of syzygy bundles.

\vspace{3mm} Let us start by fixing the notation and some basic definitions.

\begin{definition} Let $E$ be a torsion free sheaf on $\PP^N$ and set
\[\mu(E):=\frac{c_1(E)}{rk(E)}.\]
The sheaf $E$ is said to be {\em semistable} in the sense of Mumford-Takemoto if \[\mu(F)\le \mu (E)\]
for all non-zero subsheaves $F\subset E$ with
$rk(F)<rk(E)$; if strict inequality holds then $E$ is {\em
stable}.

Note that for rank $r$, torsion free sheaves $E$ on $\PP^N$, with $(c_1(E),r)=1$, the concepts of stability and semistability coincide.
\end{definition}

\vskip 2mm
\begin{notation}
Let $E$ be a rank $r$ vector bundle on $\PP^N$. We set $E_{norm}:=E(k_E)$ where  $k_E$ is the unique integer such that $c_1(E(k_E))\in \{ -r+1, \ldots , 0 \}$.
\end{notation}

\vskip 2mm For rank 2 vector bundles on $\PP^N$ we have the following useful stability criterion: a rank 2 vector bundle $E$ on $\PP^N$ is stable (resp.\ semistable) if and only if $H^0(\PP^N,E_{norm})=0$ (resp.\ $H^0(\PP^n,E_{norm}(-1))=0$).
This criterion was generalized by Hoppe in \cite{Ho}, Lemma 2.6. We have

\begin{proposition}\label{Hoppecriterio}
Let $E$ be a rank $r$ vector bundle  on $\PP^N$. The following hold:
\begin{itemize}
\item[(a)] If $H^0(X,(\wedge^qE)_{norm})=0$ for $1\leq q\leq r-1$, then $E$ is stable.
\item[(b)] $H^0(X,(\wedge^qE)_{norm}(-1))=0$ for $1\leq q\leq r-1$ if and only if  $E$ is semistable.
\end{itemize}
\end{proposition}

\vskip 2mm
\begin{remark} The conditions of Proposition~\ref{Hoppecriterio}(a) are not necessary. The simplest counterexamples are the nullcorrelation bundles $E$ on $\PP^N$ ($N$ odd) where by a nullcorrelation bundle we mean a rank $N-1$ vector bundle $E$ on $\PP^N$ ($N$ odd) defined by an exact sequence
\[
\xymatrix{
0\ar[r] &\cO_{\PP^N}(-1)\ar[r] &\Omega^1_{\PP^N}(1)\ar[r]
    &E\ar[r] &0.
}
\]
$E$ is a stable vector bundle of rank $N-1$ on $\PP^N$ ($N$ odd) and $H^0(\PP^N,(\wedge^2E)_{norm})\ne 0$ (in fact, $(\wedge ^2E)_{norm}$ contains $\cO_{\PP^N}$ as a direct summand).
\end{remark}

\begin{definition} A \emph{syzygy bundle} $E_{d_1,d_2,\ldots,d_n}$ on $\PP^N$ is a vector bundle defined as the kernel of an epimorphism
\[
\xymatrix{
\phi=(f_1,\ldots ,f_n):\oplus _{i=1}^n\cO _{\PP^N}(-d_{i})
    \ar[r] &\cO_{\PP^N},
}
\]
where $(f_1,\ldots,f_n)\subset K[X_0,X_1,\ldots,X_N]$ is an \mbox{$\frak{m}$-primary} ideal, and $d_i=deg(f_i)$.

When $d_1=d_2=\cdots =d_n=d$, we write $E_{d,n}$ instead of $E_{d_1,\ldots ,d_n}$.
\end{definition}

Let $E_{d_1,d_2,\ldots ,d_n}$ be a syzygy bundle on $\PP^N$. Since $(f_1,\ldots ,f_n)$ is an \mbox{$\frak{m}$-primary} ideal, we have $n\ge N+1$. Note also that $rank(E_{d_1,d_2,\ldots ,d_n})=n-1$, $c_1(E_{d_1,d_2,\ldots ,d_n})=-\sum_{i=1}^nd_i$ and the slope of $E_{d_1,d_2,\ldots ,d_n}$ is
\[\mu (E_{d_1,d_2,\ldots,d_n})=-\frac{\sum_{i=1}^nd_i}{n-1}.\]

\begin{definition}
A \emph{syzygy sheaf} $E_{d_1,\ldots ,d_n}$ on $\PP^N$ is a coherent sheaf defined as the kernel of  a morphism
\[
\xymatrix{
\phi:\oplus_{i=1}^n\cO_{\PP^N}(-d_{i})
    \ar[rr]^-{f_1,\ldots,f_n}
    &&\cO_{\PP^N},
}
\]
where $f_1,\ldots ,f_n\in K[X_0,X_1,\ldots ,X_N]$ are forms of degree $d_i=deg(f_i)$. When $d_1=d_2=\cdots=d_n=d$ we write $E_{d,n}$ instead of $E_{d_1,\ldots ,d_n}$.
\end{definition}

Let $E_{d_1,\ldots ,d_n}$ be a syzygy sheaf on $\PP^N$. By
construction, $E_{d_1,d_2,\ldots ,d_n}$ is a torsion free sheaf of rank $n-1$, locally free on $\cup_{i=1}^nD_{+}(f_i)\subset\PP^N$. Moreover, we have $c_1(E_{d_1,\ldots,d_n})=d-\sum_{i=1}^nd_i$, where $d$ is the degree of the highest common factor of $f_1,\ldots , f_n$ and hence the slope of $E_{d_1,d_2,\ldots ,d_n}$ is
\[\mu (E_{d_1,d_2,\ldots,d_n})=
    \frac{d-\sum_{i=1}^nd_i}{n-1}.\]

\vskip 2mm  In this paper, we address the following problems:

\begin{problem} Let $E_{d_1,d_2,\ldots,d_n}$ be the syzygy bundle on $\PP^N$ associated to a family $f_1,\ldots$, $f_n\in R$ of \mbox{$\mathfrak{m}$-primary} homogeneous polynomials of degree $d_i=deg(f_i)$. When is $E_{d_1,d_2,\ldots ,d_n}$ (semi)stable?
\end{problem}

\begin{problem} Let $f_1,\ldots ,f_n\in R$ be a family of \mbox{$\frak{m}$-primary} forms of the same degree $d$ and let $E_{d,n}$ the syzygy bundle associated to them. Is $E_{d,n}$ a (semi)stable vector bundle on $\PP^N$?
\end{problem}

As far as we know, there exist very few contributions to the above problems, and we summarize all of them, as well as the
techniques that have been used to prove these results.

\vskip 2mm First of all we observe that, as an easy application of Hoppe's Theorem, we obtain the following result, which also follows from \cite{BS}, Theorem 2.7.

\begin{proposition} \label{BSpindler}
Let $E_{d,N+1}$ be the syzygy bundle on $\PP^N$ associated to $N+1$ generic forms of degree $d$. Then, $E_{d,N+1}$ is stable.
\end{proposition}
\begin{proof}
Since stability is preserved by duality, it is enough to check that $F=E_{d,N+1}^{\vee}$ is stable. According to Proposition~\ref{Hoppecriterio}, it is enough to prove that
$H^0\left(\PP^N,(\wedge ^qF)_{norm}\right)=0$ for $1\le q \le N-1$. First of all, note that since $c_1(\wedge^qF) = \tbinom{N-1}{q-1}(N+1)q$, we have $(\wedge^qF)_{norm} = (\wedge^qF)(k_F) $ with $k_F<-dq$. Twisting by $\cO_{\PP^N}(k_F)$ the $q$th wedge power of the exact sequence
\[
\xymatrix{
0\ar[r] &\cO_{\PP^N}\ar[r] &\cO_{\PP^N}(d)^{N+1}\ar[r]
    &F\ar[r] &0,
}
\]
we get the long exact sequence:
\begin{multline*}
\xymatrix@C-.65em{
0\ar[r] &\cO_{\PP^N}(k_F)\ar[r]
    &\cO_{\PP^N}(k_F) \otimes \cO_{\PP^N}(d)^{N+1}\ar[r]
    &\cO_{\PP^N}(k_F) \otimes \wedge^2\left(\cO_{\PP^N}(d)^{N+1}\right)
    \ar[r] &\cdots
}
\\
\xymatrix@C-.65em{
\ar[r] &\cO_{\PP^N}(k_F) \otimes
    \wedge^{q-1}\left(\cO_{\PP^N}(d)^{N+1}\right)\ar[r]
    &\cO_{\PP^N}(k_F) \otimes \wedge^{q}\left(\cO_{\PP^N}(d)^{N+1}\right)\ar[r]
    &\wedge^{q}F(k_F)\ar[r] &0.
}
\end{multline*}
Cutting it into short exact sequences, for $2 \leq j \leq q-1$, we get:
\[
\xymatrix{
0\ar[r] &K_{q+1-i}\ar[r]
    &\wedge^i\left(\cO_{\PP^N}(d)^{N+1}\right)(k_F)\ar[r]
    &K_{q+2-i}\ar[r] &0,
}
\]
and
\[
\xymatrix{
0\ar[r] &\cO_{\PP^N}(k_F)\ar[r]
    &\cO_{\PP^N}(d+k_F)^{N+1}\ar[r]
    &K_{q-1}\ar[r] &0.
}
\]
Since line bundles on $\PP^N$ have no intermediate cohomology, taking cohomology on the above exact sequences we obtain
\[h^1(K_1)=h^2(K_2)= \cdots =
    h^{q-1}(K_{q-1})=h^q(\cO_{\PP^N}(k_F))=0, \]
where the last equality follows from the fact that $q < N$. On the other hand, since $k_F < -qd$,
\[H^0\left(\cO_{\PP^N}(k_F)
    \otimes \wedge^{q}
    \left(\cO_{\PP^N}(d)^{N+1}\right)\right)=
    H^0\left(\cO_{\PP^N}(qd+k_F)^{\tbinom{N+1}{q}}\right)=0.\]
Putting all together we get that for $1\le q \le N-1$,
\[H^0\left(\PP^N,(\wedge^qF)(k_F)\right)=
    H^0\left(\PP^N,(\wedge^qF)_{norm}\right)=0,\]
which proves that $F$, and hence $E_{d,N+1}$, is stable.
\end{proof}

\vskip 2mm Using the fact that the syzygy bundle $E_{d,\tbinom{d+N}{d}}$ on $\PP^N$ is a homogeneous bundle, to prove the stability of $E_{d,\tbinom{d+N}{d}}$ it is enough to check that the slope of any homogeneous sub-bundle of $E_{d,\tbinom{d+N}{d}}$ is less than the slope of $E_{d,\tbinom{d+N}{d}}$. In \cite{P}, the author described all the homogenous sub-bundles of  $E_{d,\tbinom{d+N}{d}}$ and she
proved

\begin{proposition} \label{homogeni} Let $E_{d,\tbinom{d+N}{d}}$ be the syzygy bundle on $\PP^N$ associated to $\tbinom{d+N}{d}$ \mbox{$K$-linearly} independent homogeneous forms of degree $d$. Then, $E_{d,\tbinom{d+N}{d}}$ is stable.
\end{proposition}
\begin{proof} See \cite{P}, Theorem 2.8.
\end{proof}

\vskip 2mm Using Klyachko results on toric bundles (\cite{K},
\cite{K1} and \cite{K2}), Brenner deduced the following nice
combinatoric criterion for the (semi)stability of  the syzygy bundle $E_{d_1,\ldots,d_n}$ in the case where the associated forms $f_1,\ldots,f_n$ are all monomials. Indeed, we have

\begin{proposition}\label{combinatoria} Let $f_i={X_0}^{i_0}{X_1}^{i_1}\cdots {X_N}^{i_N}$, $i\in I$, be a set of $\frak{m}$-primary monomials of degree $d_i=\sum_{j=0}^Ni_j$. Then the syzygy bundle $E_{d_1,\ldots,d_n}$ on $\PP^N$ associated to the $f_i$, $i\in I$, is semistable (resp. stable) if and only if for every $J \varsubsetneq I$, $|J|\ge 2$, the inequality
\begin{equation} \label{ineq1}
\frac{d_J-\sum _{i\in J}d_i}{|J|-1}\le
    \frac{-\sum _{i\in I}d_i}{|I|-1}\qquad\mbox{(resp.\ $<$)}
\end{equation}
holds, where $d_J$ is the degree of the greatest common factor of the $f_i$, $i\in J$.
\end{proposition}
\begin{proof}
See \cite{B}, Proposition 2.2 and Corollary 6.4.
\end{proof}

\begin{example}
(1) If we consider the set $I:=\big\{{X_0}^5,\,{X_1}^5,\,{X_2}^5, \,{X_0}^{2}{X_1}^{2}{X_2}\big\}$ of \mbox{$\frak{m}$-primary} monomials, inequality (\ref{ineq1}) is strictly fulfilled for any proper subset $J\varsubsetneq I$. Therefore the syzygy bundle $E$ associated to $I$ is stable.

(2) If we consider the set $I:=\big\{{X_0}^5,\,{X_1}^5,\,{X_2}^5,\,{X_0}^{4}{X_1}\big\}$ of \mbox{$\frak{m}$-primary} monomials, then for the subset
$J:=\big\{{X_0}^5,\,{X_0}^{4}{X_1} \big\}$  inequality (\ref{ineq1}) is not fulfilled. Therefore the syzygy bundle $E$ associated to $I$ is not stable. In fact, the slope of $E$ is
$\mu(E)=-20/3$ and the syzygy sheaf $F$ associated to $J$ is a
subsheaf of $E$ with slope $\mu(F)=-6$. Since $\mu(F)\nleq
\mu(E)$, we conclude that $E$ is not stable.
\end{example}
\begin{remark} \label{rem}
(a) Let $I$ be a set of $n$ \mbox{$\frak{m}$-primary} monomials of degree $d$. It easily follows from the above proposition that the syzygy bundle $E_{d,n}$ on $\PP^N$ associated $I$ is (semi)stable if and only if for every subset $J\subset I$ with $k:=|J|\ge 2$,
\begin{equation} \label{ineq3}
(d-d_J)n+d_J-dk >0 \quad (\mbox{resp.\ } \ge 0),
\end{equation}
where $d_J$ is the degree of the greatest common factor of the monomials in $J$.

(b) If we use the notation $a_{d,j}:=-\frac{jd}{j-1}$, inequality (\ref{ineq1}) can be written \[\tfrac{d_J}{k-1}+a_{d,k}\le a_{d,n}.\]
The fact that, once $d$ is fixed, the sequence $(a_{d,j})_{j\ge 2}$ is monotonically increasing  will be useful in many arguments.
\end{remark}

Due to Proposition~\ref{combinatoria}, to decide whether a syzygy bundle on $\PP^N$ associated to a set of \mbox{$\frak{m}$-primary} monomials of degree $d$ is semistable or not is a purely combinatorial problem but not yet solved, even when all monomials $f_{i}$ have the same degree. In \cite{B}, Question 7.8, Brenner asks

\begin{question}\label{Brennerquestion}
Does there exist for every $d$ and every $n\le\binom{N+d}{N}$ a family of $n$ monomials in $K[X_0,\ldots,X_N]$ of degree $d$ such that their syzygy bundle is semistable?
\end{question}

\begin{remark}
For $N=1$, $d=9$ and $n=3$ the answer to this question is negative. In fact, if we consider a family $I:=\big\{X^9,\,Y^9,\,X^\alpha Y^{9-\alpha}\big\}$, with ${\alpha\ge9-\alpha}$, i.e. ${\alpha\ge5}$, the subset $J\subset I$ with a greatest common factor of highest degree is $\big\{X^9,\,X^\alpha Y^{9-\alpha}\big\}$, its greatest common factor is $X^\alpha$, but inequality~\ref{ineq1} fails, since ${(9-\alpha)\cdot3+\alpha-9\cdot2=9-2\alpha<0}$.
\end{remark}

\section{The case $N=2$. Stability}

The goal of this section is to solve Problem \ref{pblm} and Question~\ref{Brennerquestion}, when $N=2$. As a main tool, we use the criterion given in Proposition \ref{combinatoria}. Let us sketch our strategy. Monic monomials in $K\left[X_0,X_1,X_2\right]$ of a given degree $d$ can be sketched in a triangle as in figure~\ref{trix}.
\begin{figure}[!htbp]
{\small%
\[
\xymatrix@C=\meiabasexgrande@R=\alturaxgrande@!0{
&&&&&&{X_2}^d\\
&&&&&{X_0}{X_2}^{d-1}&&{X_1}{X_2}^{d-1}\\
&&&&{X_0}^2{X_2}^{d-2}&&{X_0}{X_1}{X_2}^{d-2}&&
    {X_1}^2{X_2}^{d-2}\\
&&&*[@!60]{\cdots}&&*[@!60]{\cdots}&&
    *[@!-60]{\cdots}&&*[@!-60]{\cdots}\\
&&{{X_0}^{d-2}{X_2}^2}&&*[@!60]{\cdots}&&&&
    *[@!-60]{\cdots}&&{X_1}^{d-2}{X_2}^2\\
&{X_0}^{d-1}X_2&&{X_0}^{d-2}X_1X_2&&{\cdots}&&{\cdots}&&
    {X_0}{X_1}^{d-2}X_2&&{X_1}^{d-1}X_2\\
{X_0}^d&&{X_0}^{d-1}X_1&&{X_0}^{d-2}{X_1}^2&&
    {\cdots}&&{X_0}^2{X_1}^{d-2}&&{X_0}{X_1}^{d-1}&&{X_1}^d\\
}
\]
}%
\caption{Monic monomials in $K\left[X_0,X_1,X_2\right]$ of degree~$d$.}\label{trix}
\end{figure}
For the sake of simplicity, we can sketch the triangle in figure~\ref{trix} as shown in figure~\ref{tri}.

\begin{figure}[!htbp]
\[
\xymatrix@C=\meiabasegrande@R=\alturagrande@!0{
&&&&&&{\bullet}\\
&&&&&{\bullet}&&{\bullet}\\
&&&&{\bullet}&&{\bullet}&&
    {\bullet}\\
&&&*[@!60]{\cdots}&&*[@!60]{\cdots}&&
    *[@!-60]{\cdots}&&*[@!-60]{\cdots}\\
&&{\bullet}&&*[@!60]{\cdots}&&&&
    *[@!-60]{\cdots}&&{\bullet}\\
&{\bullet}&&{\bullet}&&{\cdots}&&{\cdots}&&
    {\bullet}&&{\bullet}\\
{\bullet}&&{\bullet}&&{\bullet}&&
    {\cdots}&&{\bullet}&&{\bullet}&&{\bullet}\\
}
\]
\caption{Simpler sketch of the monomials in figure~\ref{trix}.}\label{tri}
\end{figure}

Once arranged in this manner, the closer two monomials are, the greater is the degree of their greatest common factor.

\clearpage

\begin{proposition}\label{interval1}
For any integer $3\le n\le 18$ and  any integer $d\ge n-2$ there is a set~$I$ of~$n$  \mbox{$\frak{m}$-primary} monomials in $K[X_0,X_1,X_2]$ of degree~$d$ such that the corresponding syzygy bundle $E_{d,n}$ is stable.
\end{proposition}
\begin{proof} We apply Proposition \ref{combinatoria} and Remark \ref{rem}. So, for any integer $3\le n\le 18$ and any integer $d\ge n-2$ we explicitly give a set $I_{d,n}$ of $n$  $\frak{m}$-primary monomials in $K[X_0,X_1,X_2]$ of degree $d$ such that for every subset $J\subset I_{d,n}$ with $k:=|J|\ge 2$, we have
\begin{equation} \label{ineq4}
(d-d_J)n+d_J-dk >0,
\end{equation}
where $d_J$ is the degree of the greatest common factor of the monomials in $J$.

Let $e_0$, $e_1$ and $e_2$ be integers such that \[e_0+e_1+e_2=d,\qquad e_0\geq e_1\geq e_2\qquad\mbox{and}
    \qquad e_0-e_2\leq1.\]
In particular, ${e_0=\left\lceil\tfrac{d}{3}\right\rceil}$.

In case $n=3$, we  consider the set
$I_{d,3}:=\big\{{X_0}^d,\,{X_1}^d,\,{X_2}^d
    \big\}. $

In case $n=4$, we consider the set
$I_{d,4}:=\big\{{X_0}^d,\,{X_1}^d,\,{X_2}^d,
    \,{X_0}^{e_0}{X_1}^{e_1}{X_2}^{e_2}\big\}.$
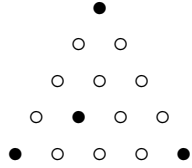
\[
\xymatrix@C=\meiabase@R=\altura@!0{
&&&&{\bullet}\\
&&&{\circ}&&{\circ}\\
&&{\circ}&&{\circ}&&{\circ}\\
&{\circ}&&{\bullet}&&{\circ}&&
    {\circ}\\
{\bullet}&&{\circ}&&{\circ}&&
    {\circ}&&{\bullet}\\
\\
&&&&I_{d,4}
}
\]

In case $n=5$, we consider the set
$I_{d,5}:=\big\{{X_0}^d,\,{X_1}^d,\,{X_2}^d,
    \,{X_0}^{e_0}{X_1}^{e_1}{X_2}^{e_2},
    \,{X_1}^{d-i}{X_2}^i\big\}$,
where ${i:=\left\lceil\tfrac{d}{2}\right\rceil}$.

In case $n=6$, we consider the set
\[I_{d,6}:=\big\{{X_0}^d,\,{X_1}^d,\,{X_2}^d,\,{X_0}^{e_0}{X_1}^{d-e_0},
    \,{X_0}^{d-e_0}{X_2}^{e_0},\,{X_1}^{e_0}{X_2}^{d-e_0}\big\}.\]

In case $n=7$, we consider the set
\[I_{d,7}:=\big\{{X_0}^d,\,{X_1}^d,\,{X_2}^d,
    \,{X_0}^{e_0}{X_1}^{e_1}{X_2}^{e_2},
    \,{X_0}^{e_0}{X_1}^{d-e_0},\,{X_0}^{d-e_0}{X_2}^{e_0},
    \,{X_1}^{e_0}{X_2}^{d-e_0}\big\}.\]

In case $n=8$, we consider the set
\begin{align*}
I_{d,8}:=\big\{&{X_0}^d,\,{X_1}^d,\,{X_2}^d,
        \,{X_0}^{e_0}{X_1}^{e_1}{X_2}^{e_2},\\
    &{X_0}^{e_0+e_1}{X_1}^{e_2},\,{X_0}^{e_2}{X_2}^{e_0+e_1},
        \,{X_1}^{e_0+e_1}{X_2}^{e_2},
        \,{X_1}^{e_0}{X_2}^{e_1+e_2}\big\}.
\end{align*}

In case $n=9$, we shall look at two cases separately: if ${d=8}$, we consider the set
\begin{align*}
I_{8,9}:=\big\{&{X_0}^8,\,{X_1}^8,\,{X_2}^8,\,{X_0}^3{X_1}^3{X_2}^2,
        \,{X_0}^6{X_1}^2,\\
    &{X_0}^2{X_2}^6,\,{X_0}^5{X_2}^3,
        \,{X_1}^6{X_2}^2,\,{X_1}^3{X_2}^5\big\};
\end{align*}
if ${d\neq8}$, let ${d=3m+t}$, with ${0\leq t<3}$, and for each ${l\in\{1,2\}}$, let ${i_l:=lm+\min(l,t)}$ and consider the set
\begin{align*}
I_{d,9}:=\big\{&{X_0}^d,\,{X_1}^d,\,{X_2}^d,
        \,{X_0}^{i_1}{X_1}^{d-i_1},
        \,{X_0}^{i_2}{X_1}^{d-i_2},\\
    &{X_0}^{d-i_1}{X_2}^{i_1},\,{X_0}^{d-i_2}{X_2}^{i_2},
        \,{X_1}^{i_1}{X_2}^{d-i_1},
        \,{X_1}^{i_2}{X_2}^{d-i_2}\big\}.
\end{align*}
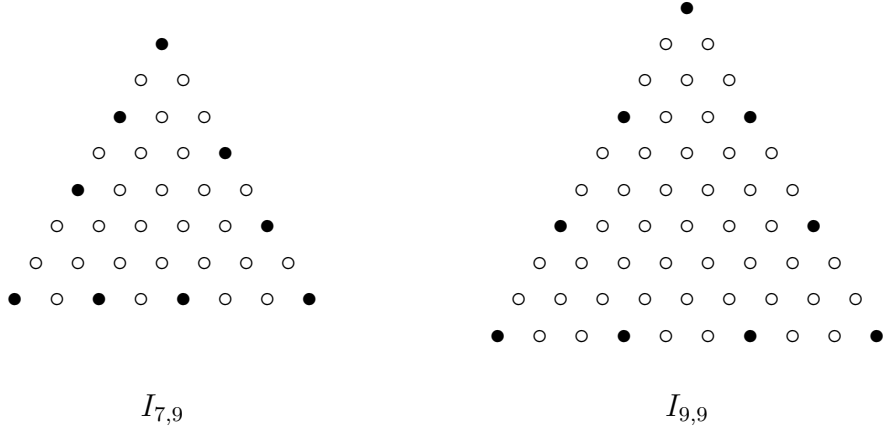
\begin{align*}
\xymatrix@C=\meiabase@R=\altura@!0{
\\
&&&&&&&{\bullet}\\
&&&&&&{\circ}&&{\circ}\\
&&&&&{\bullet}&&{\circ}&&{\circ}\\
&&&&{\circ}&&{\circ}&&{\circ}&&{\bullet}\\
&&&{\bullet}&&{\circ}&&{\circ}&&{\circ}&&{\circ}\\
&&{\circ}&&{\circ}&&{\circ}&&{\circ}&&{\circ}&&{\bullet}\\
&{\circ}&&{\circ}&&{\circ}&&{\circ}&&{\circ}&&{\circ}&&
    {\circ}\\
{\bullet}&&{\circ}&&{\bullet}&&{\circ}&&{\bullet}&&{\circ}&&
    {\circ}&&{\bullet}\\
\\
\\
&&&&&&&I_{7,9}
}
&&
\xymatrix@C=\meiabase@R=\altura@!0{
&&&&&&&&&{\bullet}\\
&&&&&&&&{\circ}&&{\circ}\\
&&&&&&&{\circ}&&{\circ}&&{\circ}\\
&&&&&&{\bullet}&&{\circ}&&{\circ}&&{\bullet}\\
&&&&&{\circ}&&{\circ}&&{\circ}&&{\circ}&&{\circ}\\
&&&&{\circ}&&{\circ}&&{\circ}&&{\circ}&&{\circ}&&{\circ}\\
&&&{\bullet}&&{\circ}&&{\circ}&&{\circ}&&{\circ}&&{\circ}&&
    {\bullet}\\
&&{\circ}&&{\circ}&&{\circ}&&{\circ}&&{\circ}&&{\circ}&&
    {\circ}&&{\circ}\\
&{\circ}&&{\circ}&&{\circ}&&{\circ}&&{\circ}&&{\circ}&&
    {\circ}&&{\circ}&&{\circ}\\
{\bullet}&&{\circ}&&{\circ}&&{\bullet}&&{\circ}&&{\circ}&&
    {\bullet}&&{\circ}&&{\circ}&&{\bullet}\\
\\
&&&&&&&&&I_{9,9}
}
\end{align*}

In case $n=10$, we shall distinguish two cases: if ${d=9}$, we consider the set
\begin{align*}
I_{9,10}:=\big\{&{X_0}^9,\,{X_1}^9,\,{X_2}^9,\,{X_0}^3{X_1}^3{X_2}^3,
        \,{X_0}^6{X_1}^3,\,{X_0}^3{X_1}^6,\\
    &{X_0}^6{X_2}^3,\,{X_0}^3{X_2}^6,\,{X_1}^6{X_2}^3,
        \,{X_1}^3{X_2}^6\big\};
\end{align*}
if $d\neq9$, let $d=5m+t$, where $0\leq t<5$, and for each $l\in\{1,2,3,4\}$, let $i_l:=lm+\min(l,t)$; consider the set
\begin{align*}
I_{d,10}:=\big\{&{X_0}^d,\,{X_1}^d,\,{X_2}^d,
        \,{X_0}^{i_2}{X_1}^{i_1}{X_2}^{d-i_1-i_2},
        \,{X_0}^{i_4}{X_1}^{d-i_4},
        \,{X_0}^{i_2}{X_1}^{d-i_2},\\
    &{X_0}^{i_3}{X_2}^{d-i_3},
        \,{X_0}^{i_1}{X_2}^{d-i_1},
        \,{X_1}^{i_2}{X_2}^{d-i_2},
        \,{X_1}^{i_4}{X_2}^{d-i_4}\big\}.
\end{align*}

In case $n=11$, we shall distinguish two cases: if $d=12$, we consider the set
\begin{align*}
I_{12,11}:=\big\{&{X_0}^{12},\,{X_1}^{12},\,{X_2}^{12},
        \,{X_0}^9{X_1}^3,\,{X_0}^6{X_1}^6,\,{X_0}^3{X_1}^9,\\
    &{X_0}^9{X_2}^3,\,{X_0}^6{X_2}^6,\,{X_0}^3{X_2}^9,
        \,{X_1}^9{X_2}^3,\,{X_1}^6{X_2}^6\big\};
\end{align*}
if $d\neq12$, let us write $d=5m+t$, where $0\leq t<5$, and for each $l\in\{1,2,3,4\}$, let $i_l:=lm+\min(l,t)$; consider the set
\begin{align*}
I_{d,11}:=\big\{&{X_0}^d,\,{X_1}^d,\,{X_2}^d,
        \,{X_0}^{i_2}{X_1}^{i_1}{X_2}^{d-i_1-i_2},
        \,{X_0}^{i_4}{X_1}^{d-i_4},
        \,{X_0}^{i_3}{X_1}^{d-i_3},\\
    &{X_0}^{i_2}{X_1}^{d-i_2},\,{X_0}^{i_3}{X_2}^{d-i_3},
        \,{X_0}^{i_1}{X_2}^{d-i_1},
        \,{X_1}^{i_2}{X_2}^{d-i_2},
        \,{X_1}^{i_4}{X_2}^{d-i_4}\big\}.
\end{align*}

In case $n=12$, we shall distinguish two cases: if $d=11$, we consider the set
\begin{align*}
I_{11,12}:=\big\{&{X_0}^{11},\,{X_1}^{11},\,{X_2}^{11},
        \,{X_0}^8{X_1}^3,\,{X_0}^8{X_2}^3,
        \,{X_0}^5{X_1}^2{X_2}^4,\,{X_0}^4{X_1}^4{X_2}^3,\\
    &{X_0}^3{X_1}^8,\,{X_0}^3{X_2}^8,\,{X_0}^2{X_1}^5{X_2}^4,\,
        \,{X_1}^8{X_2}^3,\,{X_1}^3{X_2}^8\big\};
\end{align*}
if $d\neq11$, let $d=4m+t$, where $0\leq t<4$, and for each $l\in\{1,2,3\}$, let $i_l:=lm+\min(l,t)$; consider the set
\begin{align*}
I_{d,12}:=\big\{&{X_0}^d,\,{X_1}^d,\,{X_2}^d,
        \,{X_0}^{i_3}{X_1}^{d-i_3},
        \,{X_0}^{i_2}{X_1}^{d-i_2},
        \,{X_0}^{i_1}{X_1}^{d-i_1},
        \,{X_0}^{i_3}{X_2}^{d-i_3},\\
    &{X_0}^{i_2}{X_2}^{d-i_2},\,{X_0}^{i_1}{X_2}^{d-i_1},
        \,{X_1}^{i_1}{X_2}^{d-i_1},
        \,{X_1}^{i_2}{X_2}^{d-i_2},
        \,{X_1}^{i_3}{X_2}^{d-i_3}\big\}.
\end{align*}

In case ${13\le n\le15}$, let $d=4m+t$, where $0\leq t<4$, and for each $l\in\{1,2,3\}$, let $i_l:=lm+\min(l,t)$. Consider the sets
\begin{align*}
I_{d,13}:=\big\{&{X_0}^d,\,{X_1}^d,\,{X_2}^d,
        \,{X_0}^{i_2}{X_1}^{d-i_3}{X_2}^{i_3-i_2},\\
    &{X_0}^{i_3}{X_1}^{d-i_3},\,{X_0}^{i_2}{X_1}^{d-i_2},
        \,{X_0}^{i_1}{X_1}^{d-i_1},\\
    &{X_0}^{i_3}{X_2}^{d-i_3},\,{X_0}^{i_2}{X_2}^{d-i_2},
        \,{X_0}^{i_1}{X_2}^{d-i_1},\\
    &{X_1}^{i_1}{X_2}^{d-i_1},
        \,{X_1}^{i_2}{X_2}^{d-i_2},
        \,{X_1}^{i_3}{X_2}^{d-i_3}\big\},\displaybreak[0]\\[1ex]
I_{d,14}:=\big\{&{X_0}^d,\,{X_1}^d,\,{X_2}^d,
        \,{X_0}^{i_2}{X_1}^{d-i_3}{X_2}^{i_3-i_2},
        \,{X_0}^{i_1}{X_1}^{d-i_2}{X_2}^{i_2-i_1},\\
    &{X_0}^{i_3}{X_1}^{d-i_3},\,{X_0}^{i_2}{X_1}^{d-i_2},
        \,{X_0}^{i_1}{X_1}^{d-i_1},\\
    &{X_0}^{i_3}{X_2}^{d-i_3},\,{X_0}^{i_2}{X_2}^{d-i_2},
        \,{X_0}^{i_1}{X_2}^{d-i_1},\\
    &{X_1}^{i_1}{X_2}^{d-i_1},
        \,{X_1}^{i_2}{X_2}^{d-i_2},
        \,{X_1}^{i_3}{X_2}^{d-i_3}\big\}\\
\intertext{and}
I_{d,15}:=\big\{&{X_0}^d,\,{X_1}^d,\,{X_2}^d,\\
    &{X_0}^{i_2}{X_1}^{d-i_3}{X_2}^{i_3-i_2},
        \,{X_0}^{i_1}{X_1}^{d-i_2}{X_2}^{i_2-i_1},
        \,{X_0}^{i_1}{X_1}^{d-i_3}{X_2}^{i_3-i_1},\\
    &{X_0}^{i_3}{X_1}^{d-i_3},\,{X_0}^{i_2}{X_1}^{d-i_2},
        \,{X_0}^{i_1}{X_1}^{d-i_1},\\
    &{X_0}^{i_3}{X_2}^{d-i_3},\,{X_0}^{i_2}{X_2}^{d-i_2},
        \,{X_0}^{i_1}{X_2}^{d-i_1},\\
    &{X_1}^{i_1}{X_2}^{d-i_1},
        \,{X_1}^{i_2}{X_2}^{d-i_2},
        \,{X_1}^{i_3}{X_2}^{d-i_3}\big\}.
\end{align*}

In case ${16\le n\le18}$, let $d=5m+t$, where $0\leq t<5$, and for each $l\in\{1,2,3,4\}$, let $i_l:=lm+\min(l,t)$. Consider the sets
\begin{align*}
I_{d,16}:=\big\{&{X_0}^d,\,{X_1}^d,\,{X_2}^d,
        \,{X_0}^{i_2}{X_1}^{d-i_3}{X_2}^{i_3-i_2},
        \,{X_0}^{i_4}{X_1}^{d-i_4},
        \,{X_0}^{i_3}{X_1}^{d-i_3},
        \,{X_0}^{i_2}{X_1}^{d-i_2},\\
    &{X_0}^{i_1}{X_1}^{d-i_1},
        \,{X_0}^{i_4}{X_2}^{d-i_4},
        \,{X_0}^{i_3}{X_2}^{d-i_3},
        \,{X_0}^{i_2}{X_2}^{d-i_2},
        \,{X_0}^{i_1}{X_2}^{d-i_1},\\
    &{X_1}^{i_1}{X_2}^{d-i_1},
        \,{X_1}^{i_2}{X_2}^{d-i_2},
        \,{X_1}^{i_3}{X_2}^{d-i_3},
        \,{X_1}^{i_4}{X_2}^{d-i_4}\big\},\displaybreak[0]\\[1ex]
I_{d,17}:=\big\{&{X_0}^d,\,{X_1}^d,\,{X_2}^d,
        \,{X_0}^{i_2}{X_1}^{d-i_3}{X_2}^{i_3-i_2},
        \,{X_0}^{i_2}{X_1}^{d-i_4}{X_2}^{i_4-i_2},\\
    &{X_0}^{i_4}{X_1}^{d-i_4},\,{X_0}^{i_3}{X_1}^{d-i_3},
        \,{X_0}^{i_2}{X_1}^{d-i_2},
        \,{X_0}^{i_1}{X_1}^{d-i_1},\\
    &{X_0}^{i_4}{X_2}^{d-i_4},\,{X_0}^{i_3}{X_2}^{d-i_3},
        \,{X_0}^{i_2}{X_2}^{d-i_2},
        \,{X_0}^{i_1}{X_2}^{d-i_1},\\
    &{X_1}^{i_1}{X_2}^{d-i_1},\,{X_1}^{i_2}{X_2}^{d-i_2},
        \,{X_1}^{i_3}{X_2}^{d-i_3},
        \,{X_1}^{i_4}{X_2}^{d-i_4}\big\}\\
\intertext{and}
I_{d,18}:=\big\{&{X_0}^d,\,{X_1}^d,\,{X_2}^d,
        \,{X_0}^{i_2}{X_1}^{d-i_3}{X_2}^{i_3-i_2},
        \,{X_0}^{i_2}{X_1}^{d-i_4}{X_2}^{i_4-i_2},
        \,{X_0}^{i_1}{X_1}^{d-i_3}{X_2}^{i_3-i_1},\\
    &{X_0}^{i_4}{X_1}^{d-i_4},\,{X_0}^{i_3}{X_1}^{d-i_3},
        \,{X_0}^{i_2}{X_1}^{d-i_2},
        \,{X_0}^{i_1}{X_1}^{d-i_1},\\
    &{X_0}^{i_4}{X_2}^{d-i_4},\,{X_0}^{i_3}{X_2}^{d-i_3},
        \,{X_0}^{i_2}{X_2}^{d-i_2},
        \,{X_0}^{i_1}{X_2}^{d-i_1},\\
    &{X_1}^{i_1}{X_2}^{d-i_1},\,{X_1}^{i_2}{X_2}^{d-i_2},
        \,{X_1}^{i_3}{X_2}^{d-i_3},
        \,{X_1}^{i_4}{X_2}^{d-i_4}\big\}.
\end{align*}
\[
\xymatrix@C=\meiabase@R=\altura@!0{
&&&&&&&&&&&&&&&&{\bullet}\\
&&&&&&&&&&&&&&&{\circ}&&{\circ}\\
&&&&&&&&&&&&&&{\circ}&&{\circ}&&{\circ}\\
&&&&&&&&&&&&&{\circ}&&{\circ}&&{\circ}&&{\circ}\\
&&&&&&&&&&&&{\bullet}&&{\circ}&&{\circ}&&{\circ}&&{\bullet}\\
&&&&&&&&&&&{\circ}&&{\circ}&&{\circ}&&{\circ}&&{\circ}&&
    {\circ}\\
&&&&&&&&&&{\circ}&&{\circ}&&{\circ}&&{\circ}&&{\circ}&&
    {\circ}&&{\circ}\\
&&&&&&&&&{\bullet}&&{\circ}&&{\circ}&&{\circ}&&{\circ}&&
    {\circ}&&{\circ}&&{\bullet}\\
&&&&&&&&{\circ}&&{\circ}&&{\circ}&&{\circ}&&{\circ}&&
    {\circ}&&{\circ}&&{\circ}&&{\circ}\\
&&&&&&&{\circ}&&{\circ}&&{\circ}&&{\circ}&&{\circ}&&
    {\circ}&&{\circ}&&{\circ}&&{\circ}&&{\circ}\\
&&&&&&{\bullet}&&{\circ}&&{\circ}&&{\bullet}&&{\circ}&&
    {\circ}&&{\bullet}&&{\circ}&&{\circ}&&{\circ}&&{\bullet}\\
&&&&&{\circ}&&{\circ}&&{\circ}&&{\circ}&&{\circ}&&
    {\circ}&&{\circ}&&{\circ}&&{\circ}&&{\circ}&&{\circ}&&
    {\circ}\\
&&&&{\circ}&&{\circ}&&{\circ}&&{\circ}&&{\circ}&&
    {\circ}&&{\circ}&&{\circ}&&{\circ}&&{\circ}&&{\circ}&&
    {\circ}&&{\circ}\\
&&&{\bullet}&&{\circ}&&{\circ}&&{\circ}&&{\circ}&&{\circ}&&
    {\bullet}&&{\circ}&&{\circ}&&{\circ}&&{\circ}&&
    {\circ}&&{\circ}&&{\bullet}\\
&&{\circ}&&{\circ}&&{\circ}&&{\circ}&&{\circ}&&{\circ}&&{\circ}&&
    {\circ}&&{\circ}&&{\circ}&&{\circ}&&{\circ}&&{\circ}&&
    {\circ}&&{\circ}\\
&{\circ}&&{\circ}&&{\circ}&&{\circ}&&{\circ}&&{\circ}&&
    {\circ}&&{\circ}&&{\circ}&&{\circ}&&{\circ}&&
    {\circ}&&{\circ}&&{\circ}&&{\circ}&&{\circ}\\
{\bullet}&&{\circ}&&{\circ}&&{\bullet}&&{\circ}&&{\circ}&&
    {\bullet}&&{\circ}&&{\circ}&&{\bullet}&&{\circ}&&
    {\circ}&&{\bullet}&&{\circ}&&{\circ}&&{\circ}&&{\bullet}\\
\\
&&&&&&&&&&&&&&&&I_{16,18}
}
\]

\vskip 2mm
For any $3\le n \le 18$ and $d\ge n-2$, we consider the described set $I_{d,n}$ and for any subset $J\subset I_{d,n}$ with $k:=|J|\ge 2$, we have to check that inequality (\ref{ineq4}) is satisfied. We check the case $n=18$ and we leave the other cases to the reader.

So, assume $n=18$. In this case we use the fact that no monomial of degree~$d_J$ divides a greater number of monomials in~$I_{d,n}$ than~${X_0}^{d_J}$.

If ${0<d_J\leq i_1}$, the multiples of~${X_0}^{d_J}$ in~$I_{d,18}$ are the monomials in the set
\begin{align*}
J:=\big\{&{X_0}^d,\,{X_0}^{i_2}{X_1}^{d-i_3}{X_2}^{i_3-i_2},
        \,{X_0}^{i_2}{X_1}^{d-i_4}{X_2}^{i_4-i_2},
        \,{X_0}^{i_1}{X_1}^{d-i_3}{X_2}^{i_3-i_1},\\
    &{X_0}^{i_4}{X_1}^{d-i_4},\,{X_0}^{i_3}{X_1}^{d-i_3},
        \,{X_0}^{i_2}{X_1}^{d-i_2},
        \,{X_0}^{i_1}{X_1}^{d-i_1},\\
    &{X_0}^{i_4}{X_2}^{d-i_4},\,{X_0}^{i_3}{X_2}^{d-i_3},
        \,{X_0}^{i_2}{X_2}^{d-i_2},
        \,{X_0}^{i_1}{X_2}^{d-i_1}\big\}.
\end{align*}
Therefore we have ${k=12}$ and
\begin{align*}
(d-d_J)n+d_J-dk&\geq18(d-d_J)+d_J-12d=6d-17d_J\geq6d-17i_1\geq\\
    &\geq13m+6t-17\min(1,t)\geq13m-11>0.
\end{align*}

If ${i_1<d_J\leq i_2}$, the multiples of~${X_0}^{d_J}$ in~$I_{d,18}$ are the monomials in the set
\begin{align*}
J:=\big\{&{X_0}^d,\,{X_0}^{i_2}{X_1}^{d-i_3}{X_2}^{i_3-i_2},
        \,{X_0}^{i_2}{X_1}^{d-i_4}{X_2}^{i_4-i_2},\\
    &{X_0}^{i_4}{X_1}^{d-i_4},\,{X_0}^{i_3}{X_1}^{d-i_3},
        \,{X_0}^{i_2}{X_1}^{d-i_2},\\
    &{X_0}^{i_4}{X_2}^{d-i_4},\,{X_0}^{i_3}{X_2}^{d-i_3},
        \,{X_0}^{i_2}{X_2}^{d-i_2}\big\}.
\end{align*}
Therefore we have ${k=9}$ and
\begin{align*}
(d-d_J)n+d_J-dk&=18(d-d_J)+d_J-9d=9d-17d_J\geq9d-17i_2\geq\\
    &\geq11m+9t-17\min(2,t)\geq11m-16>0.
\end{align*}

If ${i_2<d_J\leq i_3}$, the multiples of~${X_0}^{d_J}$ in~$I_{d,18}$ are the monomials in the set
\begin{align*}
J:=\big\{&{X_0}^d,\,{X_0}^{i_4}{X_1}^{d-i_4},
        \,{X_0}^{i_3}{X_1}^{d-i_3},\,{X_0}^{i_4}{X_2}^{d-i_4},
        \,{X_0}^{i_3}{X_2}^{d-i_3}\big\}.
\end{align*}
Therefore we have ${k=5}$ and
\begin{align*}
(d-d_J)n+d_J-dk&=18(d-d_J)+d_J-5d=13d-17d_J\geq13d-17i_3\geq\\
    &\geq14m+13t-17\min(3,t)\geq14m-12>0.
\end{align*}

If ${i_3<d_J\leq i_4}$, the multiples of~${X_0}^{d_J}$ in~$I_{d,18}$ are the monomials in the set
\begin{align*}
J:=\big\{&{X_0}^d,\,{X_0}^{i_4}{X_1}^{d-i_4},
        \,{X_0}^{i_4}{X_2}^{d-i_4}\big\}.
\end{align*}
Therefore we have ${k=3}$ and
\begin{align*}
(d-d_J)n+d_J-dk&=18(d-d_J)+d_J-3d=15d-17d_J\geq15d-17i_4\geq\\
    &\geq7m-2t\geq7m-8>0.
\end{align*}

If ${i_4<d_J<d}$, the only multiple of~${X_0}^{d_J}$ in~$I_{d,18}$ is ${X_0}^d$ and we have nothing to check.

Thus, we conclude that the stability is guaranteed in all the
cases.
\end{proof}

\begin{proposition}\label{interval2}
For any integers $n$ and $d$ such that  $18 < n\le d+2$, there is a set $I_{d,n}$ of $n$  $\frak{m}$-primary monomials in $K[X_0,X_1,X_2]$ of degree $d$ such that the corresponding syzygy bundle $E_{d,n}$ is stable.
\end{proposition}
\begin{proof}
For each integer $j\ge 1$, let $T_j:=\tbinom{j+1}{2}$ be the
$j$th triangular number. Choose $j$ such that ${T_{j+2}\leq n<T_{j+3}}$, and write $n=T_{j+2}+r$, with $0\le r\le j+2$. Since $n>18$, we have $j\ge 3$.  Since ${n\leq d+2}$, we get ${T_{j+2}\leq d+2}$, and therefore ${2d-j^2-5j-2\geq0}$.

\vskip 2mm

From now until the end of this proof we shall adopt the following strategy:

\vskip 2mm \noindent {\bf Strategy:} \label{X_0}
For each given~$d$ and~$n$, we choose a set of~$n$ monomials~$I_{d,n}$ such that for ${0<d_J<d}$, no monomial of degree~$d_J$ divides a greater number of monomials in~$I_{d,n}$ than~${X_0}^{d_J}$.

We write ${d=m(j+1)+t}$, where ${0\leq t<j+1}$. Note that, since ${2d\geq j^2+5j+2}$, we get ${d\geq3(j+1)+1}$, and therefore ${m\geq3}$. For each ${l\in\{1,\ldots,j\}}$, we define
\[i_l:= lm + \min(l,t). \]
We have  ${0<i_1<\cdots<i_j<d}$,
\[{d-i_j\leq i_j-i_{j-1}\leq\cdots\leq i_2-i_1\leq i_1},\]
and ${i_1-(d-i_j)\leq1}$. Set $e:=\left\lceil\tfrac{m}{2}\right\rceil$. Consider the set
\begin{align*}
I':=\big\{&{X_0}^d,\,{X_0}^{i_j}{X_1}^{d-i_j},
        \,{X_0}^{i_j}{X_2}^{d-i_j},\\
    &{X_0}^{i_{j-1}}{X_1}^{d-i_{j-1}},\,
        {X_0}^{i_{j-1}}{X_1}^{d-i_j}{X_2}^{i_j-i_{j-1}},\,
        {X_0}^{i_{j-1}}{X_2}^{d-i_{j-1}},\\
    &{X_0}^{i_{j-2}}{X_1}^{d-i_{j-2}},\,
        {X_0}^{i_{j-2}}{X_1}^{d-i_{j-1}}
        {X_2}^{i_{j-1}-i_{j-2}},\, {X_0}^{i_{j-2}}{X_1}^{d-i_j}{X_2}^{i_j-i_{j-2}},\,
        {X_0}^{i_{j-2}}{X_2}^{d-i_{j-2}},\\
    &\ldots\\
    &{X_0}^{i_1}{X_1}^{d-i_1},\,
        {X_0}^{i_1}{X_1}^{d-i_2}{X_2}^{i_2-i_1},\\
    &\qquad {X_0}^{i_1}{X_1}^{d-i_3}{X_2}^{i_3-i_1},
        \ldots,{X_0}^{i_1}{X_1}^{d-i_j}{X_2}^{i_j-i_1},\,
        {X_0}^{i_1}{X_2}^{d-i_1},\\
    &{X_1}^d,\,{X_1}^{i_j}{X_2}^{d-i_j},\ldots,
        {X_1}^{i_1}{X_2}^{d-i_1},\,{X_2}^d\big\},
\end{align*}
and the sequence
\begin{align*}
\big(&{X_0}^{i_j+e}{X_1}^{d-i_j-e},\,
        {X_0}^e{X_2}^{d-e},\,
        {X_1}^{i_j+e}{X_2}^{d-i_j-e},\\
    &{X_0}^{i_{j-1}+e}{X_1}^{d-i_{j-1}-e},\,
        {X_0}^{i_1+e}{X_2}^{d-i_1-e},\,
        {X_1}^{i_{j-1}+e}{X_2}^{d-i_{j-1}-e},\\
    &\ldots,\\
    &{X_0}^{i_{j-q}+e}{X_1}^{d-i_{j-q}-e},\,
        {X_0}^{i_q+e}{X_2}^{d-i_q-e},\,
        {X_1}^{i_{j-q}+e}{X_2}^{d-i_{j-q}-e}\big),
\end{align*}
where ${q:=\left\lceil\tfrac{j-1}{3}\right\rceil}$. Let $I''$ be the set of the first~$r$ monomials in this sequence, and let ${I_{d,n}=I'\cup I''}$. Since $I'$ has $T_{j+2}$ monomials, the number of monomials in~$I_{d,n}$ is~$n$.

\vskip 2mm

For ${1\leq l\leq j}$, let $J_l$ be the set of monomials in~$I'$ that are multiples of~${X_0}^{i_l}$. We have
\begin{align*}
J_l:=\big\{&{X_0}^d,\,{X_0}^{i_j}{X_1}^{d-i_j},\,{X_0}^{i_j}{X_2}^{d-i_j},\\
    &{X_0}^{i_{j-1}}{X_1}^{d-i_{j-1}},\,
        {X_0}^{i_{j-1}}{X_1}^{d-i_j}{X_2}^{i_j-i_{j-1}},\,
        {X_0}^{i_{j-1}}{X_2}^{d-i_{j-1}},\\
    &{X_0}^{i_{j-2}}{X_1}^{d-i_{j-2}},\,
        {X_0}^{i_{j-2}}{X_1}^{d-i_{j-1}}{X_2}^{i_{j-1}-i_{j-2}},\, {X_0}^{i_{j-2}}{X_1}^{d-i_j}{X_2}^{i_j-i_{j-2}},\,
        {X_0}^{i_{j-2}}{X_2}^{d-i_{j-2}},\\
    &\ldots\\
    &{X_0}^{i_l}{X_1}^{d-i_l},\,
        {X_0}^{i_l}{X_1}^{d-i_{l+1}}{X_2}^{i_{l+1}-i_l},\\
    &\qquad {X_0}^{i_l}{X_1}^{i_{l+2}}
        {X_2}^{i_{l+2}-i_l},
        \ldots,{X_0}^{i_l}{X_1}^{d-i_j}{X_2}^{i_j-i_l},\,
        {X_0}^{i_l}{X_2}^{d-i_l}\big\}
\end{align*}
and $\left|J_l\right|=T_{j-l+2}$.

\vskip 2mm

We  distinguish two cases.

\noindent {\bf Case 1:} ${n=T_{j+2}}$. Since we are following the strategy mentioned above, for ${0<d_J<d}$, we only have to check inequality~(\ref{ineq4}) for multiples of~${X_0}^{d_J}$.

If ${0<d_J\leq i_1}$, the multiples of~${X_0}^{d_J}$ in~$I_{d,n}$ are among the monomials in the set $J_1$.
Therefore if~$k$ is the number of multiples of~${X_0}^{d_J}$, we have ${k=T_{j+1}}$ and
\begin{align*}
(d-d_J)n+d_J-dk&=(d-d_J)T_{j+2}+d_J-dT_{j+1}=\\
    &=d(j+2)-d_JT_{j+2}+d_J\geq\\
    &\geq (m(j+1)+t)(j+2)-i_1T_{j+2}+i_1 =\\
    &= (m(j+1)+t)(j+2)-(m+\min(1,t))T_{j+2}+i_1.
\end{align*}
This last expression takes the following forms, depending on the different values of~$t$:
\begin{align*}
&\mbox{- for $t=0$, }&&\tfrac{m}{2}j(j+1);\\
&\mbox{- for $t>0$, }&&\tfrac{(j+2)}{2}
        \big(2t+(m-1)(j-1)-4 \big)+m+1.
\end{align*}
These expressions are positive in both cases because $j\ge 3$ and $m\ge 3$. So inequality (\ref{ineq4}) is strictly satisfied.

If ${i_l<d_J\leq i_{l+1}}$, for ${1\leq l\leq j-1}$, the multiples of~${X_0}^{d_J}$ in~$I_{d,n}$ are the monomials in the set $J_{l+1}$. Therefore we have ${k=T_{j+1-l}}$ and
\begin{align*}
(d-d_J)n+d_J-dk&=(d-d_J)T_{j+2}+d_J-dT_{j+1-l}=\\
    &=d(T_{j+2}-T_{j+1-l})-d_J(T_{j+2})+d_J \geq\\
    &\geq d(T_{j+2}-T_{j+1-l})-i_{l+1}(T_{j+2})+i_{l+1}.
\end{align*}
This last expression takes the following forms, depending on the different values of~$t$:
\begin{align*}
\mbox{- for }&t\leq l+1,\\
    &\tfrac{1}{2}(m-1)lj(j-l) +\tfrac{1}{2}(m-1)j(j-l)
        + \tfrac{1}{2}(m-2)l(j-l) +\\
    &\qquad +  \tfrac{1}{2}(m-2)(j-l)
    + \tfrac{1}{2}(l-1)^2(j-l) + (l-1)(j-l)\\
    &\qquad + \tfrac{1}{2}\left((j-l)^2 +3(j-l)\right)
        (l+1-t)>0;\displaybreak[0]\\
\mbox{- for }&t>l+1,\\
    &\tfrac{1}{2}(m-1)lj(j-l) + \tfrac{1}{2}(m-1)j(j-l) +
        \tfrac{1}{2}(m-2)l(j-l) +\\
    &\qquad +  \tfrac{1}{2}(m-3)(j-l)+ \tfrac{1}{2}l^2(j-l)+\\
    &\qquad +\tfrac{1}{2}\left(2l(j-l) +l^2 +2j +3l +4 \right)
        (t-l-1)>0.
\end{align*}
Therefore inequality (\ref{ineq4}) is strictly satisfied.

If ${i_j<d_J<d}$, the only multiple of~${X_0}^{d_J}$ in~$I_{d,n}$ is ${X_0}^d$ and there is nothing to check.

Therefore all possible values of $d_J$ are verified, and hence the syzygy bundle~$E_{d,n}$ is stable.

\vskip 2mm \noindent {\bf Case 2:} ${n>T_{j+2}}$. 

Here is a picture of~$I_{d,n}$ in case ${n=19}$ and ${d=20}$. In this case, we get ${j=3}$ and ${d=5(j+1)}$, therefore ${m=5}$, ${t=0}$ and ${e=3}$.

\[
\xymatrix@C=\meiabase@R=\altura@!0{
&&&&&&&&&&&&&&&&&&&&{\bullet}\\
&&&&&&&&&&&&&&&&&&&{\circ}&&{\circ}\\
&&&&&&&&&&&&&&&&&&{\circ}&&{\circ}&&{\circ}\\
&&&&&&&&&&&&&&&&&{\bullet}&&{\circ}&&{\circ}&&{\circ}\\
&&&&&&&&&&&&&&&&{\circ}&&{\circ}&&{\circ}&&{\circ}&&{\circ}\\
&&&&&&&&&&&&&&&{\bullet}&&{\circ}&&{\circ}&&{\circ}&&{\circ}&&
    {\bullet}\\
&&&&&&&&&&&&&&{\circ}&&{\circ}&&{\circ}&&{\circ}&&{\circ}&&
    {\circ}&&{\circ}\\
&&&&&&&&&&&&&{\circ}&&{\circ}&&{\circ}&&{\circ}&&{\circ}&&
    {\circ}&&{\circ}&&{\circ}\\
&&&&&&&&&&&&{\circ}&&{\circ}&&{\circ}&&{\circ}&&{\circ}&&
    {\circ}&&{\circ}&&{\circ}&&{\circ}\\
&&&&&&&&&&&{\circ}&&{\circ}&&{\circ}&&{\circ}&&{\circ}&&
    {\circ}&&{\circ}&&{\circ}&&{\circ}&&{\circ}\\
&&&&&&&&&&{\bullet}&&{\circ}&&{\circ}&&{\circ}&&{\circ}&&
    {\bullet}&&{\circ}&&{\circ}&&{\circ}&&{\circ}&&{\bullet}\\
&&&&&&&&&{\circ}&&{\circ}&&{\circ}&&{\circ}&&{\circ}&&
    {\circ}&&{\circ}&&{\circ}&&{\circ}&&{\circ}&&{\circ}&&
    {\circ}\\
&&&&&&&&{\circ}&&{\circ}&&{\circ}&&{\circ}&&{\circ}&&{\circ}&&
    {\circ}&&{\circ}&&{\circ}&&{\circ}&&{\circ}&&{\circ}&&
    {\circ}\\
&&&&&&&{\circ}&&{\circ}&&{\circ}&&{\circ}&&{\circ}&&{\circ}&&
    {\circ}&&{\circ}&&{\circ}&&{\circ}&&{\circ}&&{\circ}&&
    {\circ}&&{\circ}\\
&&&&&&{\circ}&&{\circ}&&{\circ}&&{\circ}&&{\circ}&&{\circ}&&
    {\circ}&&{\circ}&&{\circ}&&{\circ}&&{\circ}&&{\circ}&&
    {\circ}&&{\circ}&&{\circ}\\
&&&&&{\bullet}&&{\circ}&&{\circ}&&{\circ}&&{\circ}&&
    {\bullet}&&{\circ}&&{\circ}&&{\circ}&&{\circ}&&{\bullet}&&
    {\circ}&&{\circ}&&{\circ}&&{\circ}&&{\bullet}\\
&&&&{\circ}&&{\circ}&&{\circ}&&{\circ}&&{\circ}&&{\circ}&&
    {\circ}&&{\circ}&&{\circ}&&{\circ}&&{\circ}&&{\circ}&&
    {\circ}&&{\circ}&&{\circ}&&{\circ}&&{\circ}\\
&&&{\circ}&&{\circ}&&{\circ}&&{\circ}&&{\circ}&&{\circ}&&
    {\circ}&&{\circ}&&{\circ}&&{\circ}&&{\circ}&&{\circ}&&
    {\circ}&&{\circ}&&{\circ}&&{\circ}&&{\circ}&&{\circ}\\
&&{\circ}&&{\circ}&&{\circ}&&{\circ}&&{\circ}&&{\circ}&&
    {\circ}&&{\circ}&&{\circ}&&{\circ}&&{\circ}&&{\circ}&&
    {\circ}&&{\circ}&&{\circ}&&{\circ}&&{\circ}&&{\circ}&&
    {\bullet}\\
&{\circ}&&{\circ}&&{\circ}&&{\circ}&&{\circ}&&{\circ}&&
    {\circ}&&{\circ}&&{\circ}&&{\circ}&&{\circ}&&{\circ}&&
    {\circ}&&{\circ}&&{\circ}&&{\circ}&&{\circ}&&{\circ}&&
    {\circ}&&{\circ}\\
{\bullet}&&{\circ}&&{\bullet}&&{\circ}&&{\circ}&&{\bullet}&&
    {\circ}&&{\bullet}&&{\circ}&&{\circ}&&{\bullet}&&{\circ}&&
    {\circ}&&{\circ}&&{\circ}&&{\bullet}&&{\circ}&&{\circ}&&
    {\circ}&&{\circ}&&{\bullet}\\
\\
&&&&&&&&&&&&&&&&&&&&I_{20,19}
}
\]

Let ${n=T_{j+2}+r}$, with ${0< r\leq j+2}$. We have ${d+2\geq T_{j+2}+1}$. From here, if ${j>3}$, we get
\[2d\geq j^2+5j+4\geq9j+4\geq8(j+1).\]
In case ${j=3}$, since ${d\geq17}$, we have ${d\geq4(j+1)+1}$. In any case, ${m\geq4}$.

\vskip 2mm

We  distinguish three subcases.

\vskip 2mm \noindent {\bf Case 2.1:} $r=3s+1$, with ${s\geq0}$.

If ${0<d_J\leq e}$, the multiples of~${X_0}^{d_J}$ in~$I_{d,n}$ are the monomials in the set
\begin{align*}
J:=J_1\cup\big\{
    &{X_0}^{i_j+e}{X_1}^{d-i_j-e},\ldots,
        {X_0}^{i_{j-s}+e}{X_1}^{d-i_{j-s}-e},\\
    &{X_0}^e{X_2}^{d-e},\ldots,
        {X_0}^{i_{s-1}+e}{X_2}^{d-i_{s-1}-e}\big\}.
\end{align*}
Therefore if~$k$ is the number of multiples of~${X_0}^{d_J}$, we have ${k=T_{j+1}+2s+1}$, and
\begin{align*}
(d-d_J)n+d_J-dk&=(d-d_J)(T_{j+2}+3s+1)+d_J-d(T_{j+1}+2s+1)=\\
    &=d(j+2+s)-d_J(T_{j+2}+3s)\geq\\
    &\geq (m(j+1)+t)(j+2+s)-e(T_{j+2}+3s)\geq\\
    &\geq (m(j+1)+t)(j+2+s)-\tfrac{m+1}{2}(T_{j+2}+3s)=\\
    &=\tfrac{1}{4}(3m-1)j^2+\tfrac{1}{4}(7m-5)(j-2)+4(m-1)+\\
        &\qquad+\tfrac{1}{2}(2m(j-2)+3(m-1))s+t(j+2+s)>0.
\end{align*}

If ${e<d_J\leq i_1}$, the multiples of~${X_0}^{d_J}$ in~$I_{d,n}$ are the monomials in the set
\begin{align*}
J:=J_1\cup\big\{
    &{X_0}^{i_j+e}{X_1}^{d-i_j-e},\ldots,
        {X_0}^{i_{j-s}+e}{X_1}^{d-i_{j-s}-e},\\
    &{X_0}^{i_1+e}{X_2}^{d-i_1-e},\ldots,
        {X_0}^{i_{s-1}+e}{X_2}^{d-i_{s-1}-e}\big\}.
\end{align*}
Therefore we have ${k=T_{j+1}+\max(2s,1)}$, and
\begin{align*}
(d-d_J)n+d_J-dk&=(d-d_J)(T_{j+2}+3s+1)+d_J-d(T_{j+1}+
        \max(2s,1))=\\
    &=d(j+3+3s-\max(2s,1))-d_J(T_{j+2}+3s)\geq\\
    &\geq d(j+3+3s-\max(2s,1))-i_1(T_{j+2}+3s)=\\
    &=(m(j+1)+t)(j+3+3s-\max(2s,1))-\\
        &\qquad-(m+\min(1,t))(T_{j+2}+3s).
\end{align*}
This last expression takes the following forms, depending on the different values of~$s$ and~$t$:
\begin{align*}
&\mbox{- for $s=t=0$, }&&\tfrac{m}{2}(j+2)(j-1);\\
&\mbox{- for $s=0$ and $t>0$, }&&\tfrac{1}{2}
        \big((m-1)(j-2)^2+(5m-7)(j-2)+4(m-3)\big)+\\
    &&&\qquad +(t-1)(j+2);\\
&\mbox{- for $s>0$ and $t=0$, }&&m\left(\tfrac{1}{2}
        j(j+3)+(j-2)s\right);\\
&\mbox{- for $s>0$ and $t>0$, }&&\tfrac{1}{2}
        \big((m-1)j^2+(3m-5)(j-2)\big)+3(m-3)+\\
    &&&\qquad +(j+3+s)(t-1)+m(j-2)s+j+4-2s.
\end{align*}
These expressions are positive in all cases because $m\ge 3$,
$j\ge 3$ and  $ s \le \frac{j+1}{3}$. So, inequality~(\ref{ineq3}) is strictly satisfied.

If ${i_l<d_J\leq i_l+e}$, for ${1\leq l\leq j-2}$, the multiples of~${X_0}^{d_J}$ in~$I_{d,n}$ are the monomials in the set
\begin{align*}
J:=J_{l+1}\cup\big\{
    &{X_0}^{i_j+e}{X_1}^{d-i_j-e},\ldots,
        {X_0}^{i_a+e}{X_1}^{d-i_a-e},\\
    &{X_0}^{i_l+e}{X_2}^{d-i_l-e},\ldots,
        {X_0}^{i_{s-1}+e}{X_2}^{d-i_{s-1}-e}\big\},
\end{align*}
where ${a=\max(j-s,l)}$ and the second line is understood to be empty if ${s\leq l}$. Therefore we have ${k=T_{j+1-l}+\min(s+1,j+1-l)+\max(s-l,0)}$, and
\begin{align*}
(d-d_J)n+d_J-dk&=(d-d_J)(T_{j+2}+3s+1)+d_J-\\
    &\qquad -d\big(T_{j+1-l}+1+\min(s,j-l)+\max(s-l,0)\big)=\\
    &=d\big(T_{j+2}-T_{j+1-l}+3s
        -\min(s,j-l)-\max(s-l,0)\big)-\\
    &\qquad -d_J(T_{j+2}+3s)\geq\\
    &\geq d\big(T_{j+2}-T_{j+1-l}+3s-\min(s,j-l)
        -\max(s-l,0)\big)-\\
    &\qquad -\left(i_l+\tfrac{m+1}{2}\right)(T_{j+2}+3s).
\end{align*}
This last expression takes the following forms, depending on the different values of~$j$, $l$, $s$ and~$t$:
\begin{align*}
\mbox{- for }&s\leq j-l,\ s\leq l,\ t\leq l,\\
    &\tfrac{1}{2}(m-2)lj(j-l-1) + \tfrac{3}{4}(m-1)j(j-l) +
        \tfrac{1}{4}(m-4)l(j-l) +\\
    &\qquad +\tfrac{3}{4}(m-4)l^2+\tfrac{7}{4}(m-3)(j-l)+
        \tfrac{1}{4}(m-4)l+\tfrac{1}{2}m+\\
    &\qquad +\tfrac{1}{2}(l-1)(j-l-2)^2 +(l-1)^2(j-l-2)
        +j(j-l-2) +\\
    &\qquad +\tfrac{3}{4}l(j-l-2)+
        \tfrac{15}{4}l(l-1)+3(j-l) +\tfrac{5}{2} +\\
    &\qquad + \tfrac{1}{2}(4m(j-l)+m-3)s +l(m+1)(j-l-s)\\
    &\qquad + \tfrac{1}{2}(l-t)\big((j-l)^2+3(j-l)+2+2s\big)
        >0;\displaybreak[0]\\
\mbox{- for }&s\leq j-l,\ s\leq l,\ t>l,\\
    &\tfrac{1}{2}(m-4)jl(j-l-1) +\tfrac{3}{4}(m-1)j(j-l)+
        \tfrac{1}{4}(m-4)lj +\tfrac{1}{2}ml^2+\\
    &\qquad +\tfrac{7}{4}(m-1)(j-l)+\tfrac{1}{4}ml+
        \tfrac{1}{2}(m-4)+\tfrac{3}{2}lj(j-l-2)+\\
    &\qquad +\tfrac{1}{2}(l-1)^2(j-l)+
        \tfrac{1}{2}j(j-l-2)+\tfrac{1}{4}l(j-l-2)+j+\\
    &\qquad  +\tfrac{7}{4}l(l-1)+
        \tfrac{1}{2}+
    \tfrac{1}{2}(t-l)\big((2j-l)l+2j+3l+4+4s\big)+\\
    &\qquad +\tfrac{1}{2}(4m(j-l)+m-3)s+
        (m+1)l(j-l-s)>0;\displaybreak[0]\\
\mbox{- for }&s\leq j-l, s>l, t\leq l,\\
    &\tfrac{1}{2}(m-1)lj(j-l) +\tfrac{3}{4}(m-1)j(j-l)+
        \tfrac{7}{4}(m-3)l(j-l) +\tfrac{1}{4}ml^2+\\
    &\qquad +\tfrac{1}{4}(m-2)(3j+4l+2) +
        \tfrac{1}{2}l^2(j-l)+\tfrac{1}{2}j^2 +l(j-l) +\\
    &\qquad +\tfrac{1}{4}l(l-1)+\tfrac{1}{4}(j-l-2) +
        \tfrac{1}{2}(2m(j-2l)+m-3)(s-l)+\\
    &\qquad + \big((m+2)l+m\big)(j-l-s) +\\
    &\qquad + \tfrac{1}{2}\left((j-l)^2+3(j-2l)+l+2+4s\right)
        (l-t)>0;\displaybreak[0]\\
\mbox{- for }&s\leq j-l,\ s>l,\ t>l,\\
    &\tfrac{1}{2}(m-1)lj(j-l) +\tfrac{3}{4}(m-1)j(j-l)+
        \tfrac{7}{4}(m-3)l(j-l) +\tfrac{1}{4}ml^2 +\\
    &\qquad +\tfrac{1}{4}(m-1)(3j+4l+2) +
        \tfrac{1}{2}l^2(j-l)+\tfrac{1}{2}j(j-2l-1) +\\
    &\qquad +\tfrac{1}{2}(2m(j-2l)+m-3)(s-l)++\tfrac{7}{4}(l-1)+
        \tfrac{3}{4}\\
    &\qquad +\big((m+2)l+m\big)(j-l-s)+2l(j-2l)+\tfrac{13}{4}l(l-1) +\\
    &\qquad +\tfrac{1}{2}\left(2(l+1)(j-l)+l^2+7l+4+2s\right)
        (t-l)>0;\displaybreak[0]\\
\mbox{- for }&s>j-l,\ s\leq l,\ t\leq l,\\
    &\tfrac{1}{2}(m-4)lj(j-l-2) +\tfrac{11}{4}(m-2)j(j-l)+
        \tfrac{5}{4}(m-1)l(2l-j) +\\
    &\qquad +\tfrac{9}{4}m(j-l) +\tfrac{1}{4}ml +
        \tfrac{1}{2}m +\tfrac{3}{2}(l-1)(j-l-2)^2  + \tfrac{1}{4}l +
        \tfrac{15}{2}+\\
    &\qquad +2(l-1)^2(j-l-2) +\tfrac{27}{4}j(j-l-2) +
        \tfrac{1}{2}l(j-l-2)+\\
    &\qquad +(l-1)^2+\tfrac{11}{4}(j-l)+
    \tfrac{1}{2}(s-j+l)( 6m(j-l) +  3m - 3 )+\\
    &\qquad + \tfrac{1}{2}(l-t)\left( (j-l)^2 +
        5(j-l) + 2 \right)>0;\displaybreak[0]\\
\mbox{- for }&s>j-l,\ s\leq l,\ t>l,\\
    &\tfrac{1}{2}(m-4)lj(j-l-2) + \tfrac{11}{4}(m-2)j(j-l) +
        \tfrac{5}{4}(m-1)l(2l-j)+ \\
    &\qquad +\tfrac{9}{4}m(j-l) +\tfrac{1}{4}ml +
        \tfrac{1}{2}m+\tfrac{3}{2}(l-1)(j-l-2)^2 +\\
    &\qquad + 2(l-1)^2(j-l-2) + \tfrac{27}{4}j(j-l-2) +
        \tfrac{1}{2}l(j-l-2) + (l-1)^2+\\
    &\qquad +\tfrac{11}{4}(j-l)+\tfrac{1}{4}l +
        \tfrac{15}{2}+\\
    &\qquad +\tfrac{1}{2}(s-j+l)(6m(j-l)+ 3m-3)+\\
    &\qquad +\tfrac{1}{2}(t-l)
        \left(2l(j-l)+l^2+5l+4+6s\right)>0.
\end{align*}
Since ${s\leq\tfrac{j+1}{3}}$, if ${l<s}$, we get ${j-l>\tfrac{2j-1}{3}\geq\tfrac{j+1}{3}\geq s}$. Therefore all possible cases are checked, and inequality~(\ref{ineq4}) is strictly satisfied.

If ${i_l+e<d_J\leq i_{l+1}}$, for ${1\leq l\leq j-2}$, the multiples of~${X_0}^{d_J}$ in~$I_{d,n}$ are the monomials in the set
\begin{align*}
J:=J_{l+1}\cup\big\{
    &{X_0}^{i_j+e}{X_1}^{d-i_j-e},\ldots,
        {X_0}^{i_a+e}{X_1}^{d-i_a-e},\\
    &{X_0}^{i_{l+1}+e}{X_2}^{d-i_{l+1}-e},\ldots,
        {X_0}^{i_{s-1}+e}{X_2}^{d-i_{s-1}-e}\big\},
\end{align*}
where ${a=\max(j-s,l+1)}$, and the second line is understood to be empty if ${s\leq l+1}$. Therefore we have ${k=T_{j+1-l}+\min(s+1,j-l)+\max(s-l-1,0)}$, and
\begin{align*}
(d-d_J)n+d_J-dk&=(d-d_J)(T_{j+2}+3s+1)+d_J-\\
    &\qquad -d\big(T_{j+1-l}+\min(s+1,j-l)+
        \max(s-l-1,0)\big)=\\
    &=d\big(T_{j+2}-T_{j+1-l}+3s+1-\min(s+1,j-l)-\\
    &\qquad -\max(s-l-1,0)\big)-d_J(T_{j+2}+3s)\geq\\
    &\geq d\big(T_{j+2}-T_{j+1-l}+3s+1-\min(s+1,j-l)-\\
    &\qquad -\max(s-l-1,0)\big)-i_{l+1}(T_{j+2}+3s).
\end{align*}
This last expression takes the following forms, depending on the different values of~$j$, $l$, $s$ and~$t$:
\begin{align*}
\mbox{- for }&s+1\leq j-l,\ s\leq l+1,\ t\leq l+1,\\
    &\tfrac{1}{2}(m-4)lj(j-l) +\tfrac{1}{2}(m-4)(j-l)(j-l-1)+
        (m-4)(j-l-1)+\\
    &\qquad +\tfrac{3}{2}lj(j-l-2)+ \tfrac{1}{2}l^2(j-l)+
        \tfrac{3}{2}(j-l)(j-l-2)+\tfrac{1}{2}l(j-l)+\\
    &\qquad + \tfrac{5}{2}(j-l-2)+ 3l^2+1+m(2(j-l-1)+1)s+\\
    &\qquad +(ml+t)(j-l-1-s)+\tfrac{1}{2}
        \left((j-l)^2+5(j-l)\right)(l+1-t)>0;\displaybreak[0]\\
\mbox{- for }&s+1\leq j-l,\ s\leq l+1,\ t>l+1,\\
    &\tfrac{1}{2}(m-4)lj(j-l) +\tfrac{1}{2}(m-4)(j-l)(j-l-1)
       + (m-4)(j-l-1)+ \\
    &\qquad  + \tfrac{3}{2}l(j-l-2)^2+ 2l(l-1)(j-l) +
        \tfrac{3}{2}(j-l)(j-l-2)+\\
    &\qquad  + \tfrac{7}{2}l(j-l-2)+ \tfrac{1}{2}(j-l-2)
        + 3(l-1)+ 2+\\
    &\qquad  +(2m(j-l-1)+m+2t)s+ (ml+3l+3)(j-l-1-s)+\\
    &\qquad + \tfrac{1}{2}\left( 2l(j-l) + l^2 + 2(j-l)
        + 5l + 4 \right)(t-l-1)>0;\displaybreak[0]\\
\mbox{- for }&s+1\leq j-l, s>l+1, t\leq l+1,\\
    &\tfrac{1}{2}(m-1)lj(j-l) +
        \tfrac{1}{2}(m-4)\left(j^2+l^2\right)+
        \tfrac{5}{2}(m-1)(j-l-1) + ml + \\
    &\qquad + \tfrac{1}{2}m+\tfrac{1}{2}(l-1)^2(j-l) +
        \tfrac{3}{2}j(j-l-1)+ 4l^2+ \tfrac{7}{2}(l-1) + 3+\\
    &\qquad + m(j-l-2)(s-l-1) + 2(ml+t)(j-l-1-s)+\\
    &\qquad + \tfrac{1}{2}\left( (j-l-2)^2
        + 11(j-l-2) + 14 \right)(l+1-t)
        + lt>0;\displaybreak[0]\\
\mbox{- for }&s+1\leq j-l,\ s>l+1,\ t>l+1,\\
    &\tfrac{1}{2}(m-4)lj(j-l) +
        \tfrac{1}{2}(m-1)\left(j^2+l^2\right)+
        \tfrac{5}{2}(m-2)(j-l-1) +ml+\\
    &\qquad +\tfrac{1}{2}m +\tfrac{3}{2}lj(j-l-2) +
        \tfrac{1}{2}l(l-1)(j-l)+ 6l^2 +\tfrac{1}{2}j +
        \tfrac{11}{2}(l-1) +\tfrac{9}{2}+\\
    &\qquad + ( m(j-l-2) + t )(s-l-1)+
        \big((2m+3)l+3\big)(j-l-1-s)+\\
    &\qquad + \tfrac{1}{2}\left( 2l(j-l) + l^2 + 2j
        + 7l + 8 \right)(t-l-1)>0;\displaybreak[0]\\
\mbox{- for }&s+1>j-l,\ s\leq l+1,\ t\leq l+1,\\
    &\tfrac{1}{2}(m-4)lj(j-l) + \tfrac{5}{2}(m-1)(j-l)^2+
        \tfrac{7}{2}(m-1)(j-l) + \tfrac{3}{2}lj(j-l-2)+\\
    &\qquad + \tfrac{1}{2}l^2(j-l) + 2(j-l)^2
        +\tfrac{1}{2}lj +\tfrac{5}{2}l^2 + j-l +\\
    &\qquad + 3m(j-l)(s-j+l-1)+\\
    &\qquad + \tfrac{1}{2}\left( (j-l)^2 + 5(j-l) \right)
        (l+1-t)>0;\displaybreak[0]\\
\mbox{- for }&s+1>j-l,\ s\leq l+1,\ t>l+1,\\
    &\tfrac{1}{2}(m-4)lj(j-l) + \tfrac{5}{2}(m-1)(j-l)^2+
        \tfrac{7}{2}(m-1)(j-l) + \\
    &\qquad + \tfrac{3}{2}lj(j-l-2)+ \tfrac{1}{2}l^2(j-l) +
        2(j-l)^2 +\tfrac{1}{2}lj+\\
    &\qquad + \tfrac{5}{2}l^2 + j-l + 3m(j-l)(s-j+l-1)+\\
    &\qquad + \tfrac{1}{2}\left( 2l(j-l) + l^2
        + 5l + 6 + 6s \right)(t-l-1)>0.
\end{align*}
Since ${s\leq\tfrac{j+1}{3}}$, if ${l+1<s}$, we get ${j-l-1>\tfrac{2j-1}{3}\geq\tfrac{j+1}{3}\geq s}$. Therefore all possible cases are checked, and inequality~(\ref{ineq4}) is strictly satisfied.

If ${i_{j-1}<d_J\leq i_{j-1}+e}$, the multiples of~${X_0}^{d_J}$ in~$I_{d,n}$ are the monomials in the set
\begin{align*}
J:=J_j\cup\big\{
    &{X_0}^{i_j+e}{X_1}^{d-i_j-e},\,
        {X_0}^{i_{j-a}+e}{X_1}^{d-i_{j-a}-e}\big\},
\end{align*}
where ${a=\min(s,1)}$. Therefore we have ${k=T_2+\min(s+1,2)}$, and
\begin{align*}
(d-d_J)n+d_J-dk&=(d-d_J)(T_{j+2}+3s+1)+d_J-d\big(T_2+\min(s+1,2)\big)=\\
    &=d\big(T_{j+2}-T_2+3s-\min(s,1)\big)-d_J(T_{j+2}+3s)\geq\\
    &\geq d\big(T_{j+2}-3+3s-\min(s,1)\big)-\\
    &\qquad -\left(i_{j-1}+\tfrac{m+1}{2}\right)(T_{j+2}+3s).
\end{align*}
This last expression takes the following forms, depending on the different values of~$j$, $s$ and~$t$:
\begin{align*}
\mbox{- for }&s\leq1,\ t\leq j-1,\\
    &\tfrac{3}{4}(m-2)j(j-1) + \tfrac{1}{2}(m-4)j +
        \tfrac{9}{2}(m-1) + \tfrac{5}{4}(j-3)^2 +
        \tfrac{11}{4}(j-3) + 4 +\\
    &\qquad +\big((m+1)(j-3)+2\big)(1-s) +
        \tfrac{1}{2}(m-3)s +
        (s+3)(j-1-t)>0;\displaybreak[0]\\
\mbox{- for }&s\leq1,\ t=j,\\
    &\tfrac{3}{4}(m-1)j(j-1) + \tfrac{1}{2}mj +
        \tfrac{9}{2}(m-1)+ (j-3)^2 + \tfrac{5}{2}j+\\
    &\qquad +(m+1)(j-3)(1-s) +
        \tfrac{1}{2}(m-3)s>0;\displaybreak[0]\\
\mbox{- for }&s>1,\ t\leq j-1,\\
    &\tfrac{3}{4}(m-2)j(j-1)+\tfrac{1}{2}(m-4)j+5(m-1)+
        \tfrac{5}{4}(j-3)^2 + \tfrac{11}{4}(j-3) + 3+\\
    &\qquad +\tfrac{1}{2}( 9(m-1) + 6 )(s-1)+
        4(j-1-t)>0;\displaybreak[0]\\
\mbox{- for }&s>1,\ t=j,\\
    &\tfrac{3}{4}(m-1)j(j-1) + \tfrac{1}{2}mj + 5(m-1)+
        (j-3)^2 + \tfrac{5}{2}(j-1) + \tfrac{3}{2}+\\
    &\qquad +\tfrac{1}{2}(9m+3)(s-1)>0.
\end{align*}
Therefore inequality~(\ref{ineq4}) is strictly satisfied.

If ${i_{j-1}+e<d_J\leq i_j}$, the multiples of~${X_0}^{d_J}$ in~$I_{d,n}$ are the monomials in the set
\[
J:=J_j\cup\big\{{X_0}^{i_j+e}{X_1}^{d-i_j-e}\big\}.
\]
Therefore we have ${k=T_2+1=4}$, and
\begin{align*}
(d-d_J)n+d_J-dk&=(d-d_J)(T_{j+2}+3s+1)+d_J-4d=\\
    &=d(T_{j+2}+3s-3)-d_J(T_{j+2}+3s)\geq\\
    &\geq d(T_{j+2}+3s-3)-i_j(T_{j+2}+3s)=\\
    &=(m(j+1)+t)(T_{j+2}+3s-3)-\\
    &\qquad -(mj+t)(T_{j+2}+3s)=\\
    &=\tfrac{1}{2}\big((m-3)j(j-1) +3j(j-3)\big) +3(j-t)
        +3ms>0.
\end{align*}

If ${i_j<d_J\leq i_j+e}$, the multiples of~${X_0}^{d_J}$ in~$I_{d,n}$ are the monomials in the set
\[
J:=\big\{{X_0}^d,\,{X_0}^{i_j+e}{X_1}^{d-i_j-e}\big\}.
\]
Therefore we have ${k=2}$, and
\begin{align*}
(d-d_J)n+d_J-dk&=(d-d_J)(T_{j+2}+3s+1)+d_J-2d=\\
    &=d(T_{j+2}+3s-1)-d_J(T_{j+2}+3s)\geq\\
    &\geq d(T_{j+2}+3s-1)-(i_j+e)(T_{j+2}+3s)\ge \\
    & \ge (m(j+1)+t)(T_{j+2}+3s-1)-\\
    &\qquad -\left(mj+t+\tfrac{m+1}{2}\right)(T_{j+2}+3s)=\\
    &=\tfrac{1}{2}\big((m-1)j^2 +3(m-1)j +3(m-1)s\big)
        +(j-t)+\\
    &\qquad +2(m-1)+1>0.
\end{align*}

If ${i_j+e<d_J<d}$, the only multiple of~${X_0}^{d_J}$ in~$I_{d,n}$ is ${X_0}^d$, and there is nothing to prove.

Therefore all possible values of $d_J$ are verified, and hence the syzygy bundle~$E_{d,n}$ is stable.

\vskip 2mm \noindent {\bf Case 2.2:} $r=3s+2$, with ${s\geq0}$. The difference between this case and the previous one is that we are adding the monomial ${X_0}^{i_s+e}{X_2}^{d-i_s-e}$ to~$I_{d,n}$. Therefore we should only worry with the cases ${0<d_J\leq i_s+e}$, since for degrees greater than ${i_s+e}$ the set $J$ of multiples of~${X_0}^{d_J}$ has the same number of elements as in the corresponding sets of the previous case, whereas the set $I_{d,n}$ has one more element. Given the fact that the sequence $\left(a_{d,j}\right)_{j\geq2}$ is monotonically increasing, inequality~(\ref{ineq4}) is strictly satisfied.

If ${0<d_J\leq e}$, the multiples of~${X_0}^{d_J}$ in~$I_{d,n}$ are among the monomials in the set
\begin{align*}
J:=J_1\cup\big\{
    &{X_0}^{i_j+e}{X_1}^{d-i_j-e},\ldots,
        {X_0}^{i_{j-s}+e}{X_1}^{d-i_{j-s}-e},\\
    &{X_0}^e{X_2}^{d-e},\ldots,
        {X_0}^{i_s+e}{X_2}^{d-i_s-e}\big\}.
\end{align*}
Therefore if~$k$ is the number of multiples of~${X_0}^{d_J}$, we have ${k=T_{j+1}+2s+2}$, and
\begin{align*}
(d-d_J)n+d_J-dk&=(d-d_J)(T_{j+2}+3s+2)+d_J-d(T_{j+1}+2s+2)=\\
    &=d(j+2+s)-d_J(T_{j+2}+3s+1)\geq\\
    &\geq (m(j+1)+t)(j+2+s)-e(T_{j+2}+3s+1)\geq\\
    &\geq (m(j+1)+t)(j+2+s)-\tfrac{m+1}{2}(T_{j+2}+3s+1)=\\
    &=\tfrac{1}{4}(3m-1)j^2+\tfrac{1}{4}(7m-5)(j-2)+
        \tfrac{7}{2}(m-2)+\tfrac{5}{2}+\\
    &\qquad+\tfrac{1}{2}(2m(j-2)+3(m-1))s+t(j+2+s)>0.
\end{align*}

If ${e<d_J\leq i_1}$, the multiples of~${X_0}^{d_J}$ in~$I_{d,n}$ are the monomials in the set
\begin{align*}
J:=J_1\cup\big\{
    &{X_0}^{i_j+e}{X_1}^{d-i_j-e},\ldots,
        {X_0}^{i_{j-s}+e}{X_1}^{d-i_{j-s}-e},\\
    &{X_0}^{i_1+e}{X_2}^{d-i_1-e},\ldots,
        {X_0}^{i_s+e}{X_2}^{d-i_s-e}\big\}.
\end{align*}
Therefore we have ${k=T_{j+1}+2s+1}$ and
\begin{align*}
(d-d_J)n+d_J-dk&=(d-d_J)(T_{j+2}+3s+2)+d_J-d(T_{j+1}+2s+1)=\\
    &=d(j+3+s)-d_J(T_{j+2}+3s+1)\geq\\
    &\geq d(j+3+s)-i_1(T_{j+2}+3s+1)=\\
    &=(m(j+1)+t)(j+3+s)-\\
        &\qquad-\big(m+\min(1,t)\big)(T_{j+2}+3s+1).
\end{align*}
This last expression takes the following forms, depending on the different values of~$t$:
\begin{align*}
&\mbox{- for $t=0$, }&&\tfrac{m}{2}(j+2)(j-1)+ m(j-2)s +mj;\\
&\mbox{- for $t>0$, }&&\tfrac{1}{2}
        \big((m-1)(j-2)^2+(7m-7)(j-2)+8(m-3)+12\big)+\\
    &&&\qquad +(m(j-3)+m-3+t)s+(t-1)(j+3).
\end{align*}
These expressions are both positive, so inequality~(\ref{ineq4}) is strictly satisfied.

If ${i_l<d_J\leq i_l+e}$, for ${1\leq l\leq s}$, we get ${j-s\geq\tfrac{2}{3}j>l}$, since ${3s+2\leq j+2}$. Therefore the multiples of~${X_0}^{d_J}$ in~$I_{d,n}$ are the monomials in the set
\begin{align*}
J:=J_{l+1}\cup\big\{
    &{X_0}^{i_j+e}{X_1}^{d-i_j-e},\ldots,
        {X_0}^{i_{j-s}+e}{X_1}^{d-i_{j-s}-e},\\
    &{X_0}^{i_l+e}{X_2}^{d-i_l-e},\ldots,
        {X_0}^{i_s+e}{X_2}^{d-i_s-e}\big\}.
\end{align*}
Therefore we have ${k=T_{j+1-l}+2s+2-l}$, and
\begin{align*}
(d-d_J)n+d_J-dk&=(d-d_J)(T_{j+2}+3s+2)+d_J-\\
    &\qquad -d(T_{j+1-l}+2s+2-l)=\\
    &=d(T_{j+2}-T_{j+1-l}+s+l)-d_J(T_{j+2}+3s+1)\geq\\
    &\geq d(T_{j+2}-T_{j+1-l}+s+l)-
        \left(i_l+\tfrac{m+1}{2}\right)(T_{j+2}+3s+1).
\end{align*}
This last expression takes the following forms, depending on the different values of~$t$:
\begin{align*}
\mbox{- for }&t\leq l,\\
    &\tfrac{1}{2}(m-1)lj(j-l) +\tfrac{3}{4}(m-1)j(j-l)+
        \tfrac{7}{4}(m-3)l(j-l) +\tfrac{1}{4}ml^2+\\
    &\qquad +\tfrac{3}{4}(m-4)j +
        \tfrac{1}{2}l^2(j-l)+\tfrac{1}{2}j(j-2l) +2l(j-2l)
        +\tfrac{13}{4}l(l-1)+\\
    &\qquad +\tfrac{7}{4}(j-l-2) +
        \tfrac{3}{2}(l+1)+\tfrac{1}{2}(2m(j-2l)+m-3)(s-l)+\\
    &\qquad + \big((m+2)l+m\big)(j-l-s) +\\
    &\qquad + \tfrac{1}{2}\left((j-l)^2+3(j-2l)+l+4+4s\right)
        (l-t)>0;\displaybreak[0]\\
\mbox{- for }&t>l,\\
    &\tfrac{1}{2}(m-1)lj(j-l) +\tfrac{3}{4}(m-1)j(j-l)+
        \tfrac{7}{4}(m-3)l(j-l) +\tfrac{1}{4}ml^2 +\\
    &\qquad +\tfrac{3}{4}(m-1)j +
        \tfrac{1}{2}l^2(j-l-2)+\tfrac{1}{2}j(j-2l-1) +\\
    &\qquad +2l(j-2l-1)+\tfrac{17}{4}l(l-1) +
        \tfrac{11}{4}(l-1)+\tfrac{3}{4}+\\
    &\qquad +\tfrac{1}{2}(2m(j-2l)+m-3)(s-l)+\\
    &\qquad +\big((m+2)l+m\big)(j-l-s)+\\
    &\qquad +\tfrac{1}{2}\left(2(l+1)(j-l)+l^2+7l+4+2s\right)
        (t-l)>0.
\end{align*}

If ${i_l+e<d_J\leq i_{l+1}}$, for ${1\leq l\leq s-1}$, the multiples of~${X_0}^{d_J}$ in~$I_{d,n}$ are the monomials in the set
\begin{align*}
J:=J_{l+1}\cup\big\{
    &{X_0}^{i_j+e}{X_1}^{d-i_j-e},\ldots,
        {X_0}^{i_{j-s}+e}{X_1}^{d-i_{j-s}-e},\\
    &{X_0}^{i_{l+1}+e}{X_2}^{d-i_{l+1}-e},\ldots,
        {X_0}^{i_s+e}{X_2}^{d-i_s-e}\big\}.
\end{align*}
Therefore we have ${k=T_{j+1-l}+2s+1-l}$, and
\begin{align*}
(d-d_J)n+d_J-dk&=(d-d_J)(T_{j+2}+3s+2)+d_J-\\
    &\qquad -d(T_{j+1-l}+2s+1-l)=\\
    &=d(T_{j+2}-T_{j+1-l}+s+1+l)-d_J(T_{j+2}+3s+1)\geq\\
    &\geq d(T_{j+2}-T_{j+1-l}+s+1+l)-i_{l+1}(T_{j+2}+3s+1).
\end{align*}
This last expression takes the following forms, depending on the different values of~$t$:
\begin{align*}
\mbox{- for }&t\leq l+1,\\
    &\tfrac{1}{2}(m-1)lj(j-l) +
        \tfrac{1}{2}(m-4)\left(j^2+l^2\right)+
        \tfrac{5}{2}(m-1)(j-l-2) + 2m + \\
    &\qquad +\tfrac{1}{2}(l-1)^2(j-l) +
        \tfrac{3}{2}j(j-l-1)+ 4l(l-1)+\tfrac{13}{2}(l-1) +
        \tfrac{5}{2}+\\
    &\qquad + m(j-l-2)(s-l-1) + 2(ml+t)(j-l-1-s)+\\
    &\qquad + \tfrac{1}{2}\left( (j-l-2)^2
        + 11(j-l-2) + 16 \right)(l+1-t)
        + lt>0;\displaybreak[0]\\
\mbox{- for }&t>l+1,\\
    &\tfrac{1}{2}(m-4)lj(j-l) +
        \tfrac{1}{2}(m-1)\left(j^2+l^2\right)+
        \tfrac{5}{2}(m-2)(j-l-2) +2m+\\
    &\qquad +\tfrac{3}{2}lj(j-l-2) +
        \tfrac{1}{2}l(l-1)(j-l)+ 6l(l-1) +\tfrac{1}{2}j +
        \tfrac{21}{2}(l-1) +\tfrac{7}{2}+\\
    &\qquad + ( m(j-l-2) + t )(s-l-1)+
        \big((2m+3)l+3\big)(j-l-1-s)+\\
    &\qquad + \tfrac{1}{2}\left( 2l(j-l) + l^2 + 2j
        + 7l + 8 \right)(t-l-1)>0.
\end{align*}
Therefore inequality~(\ref{ineq4}) is strictly satisfied.

\vskip 2mm \noindent {\bf Case 2.3:} $r=3s$, with ${s\geq1}$. The difference between this case and the previous one is that we are adding the monomial ${X_1}^{i_s+e}{X_2}^{d-i_s-e}$ to~$I_{d,n}$. Since this is no multiple of~${X_0}^{d_J}$, the set $J$ of multiples of~${X_0}^{d_J}$ has the same number of elements as in the corresponding sets of the previous case, whereas the set $I_{d,n}$ has one more element. Given the fact that the sequence $\left(a_{d,j}\right)_{j\geq2}$ is monotonically increasing (see Remark \ref{rem}(b)), inequality~(\ref{ineq4}) is strictly satisfied.

We can conclude that stability is guaranteed in all cases.
\end{proof}

\begin{proposition}\label{interval3}
For any integers $n$ and $d$ such that  $d+2< n\le 3d$ and $(n,d)\neq(5,2)$, there is a set $I_{d,n}$ of $n$ $\frak{m}$-primary monomials in $K[X_0,X_1,X_2]$ of degree $d$ such that the corresponding syzygy bundle $E_{d,n}$ is stable. For $(n,d)=(5,2)$, there are $5$  $\frak{m}$-primary monomials in $K[X_0,X_1,X_2]$ of degree $2$ such that the corresponding syzygy bundle $E_{2,5}$ is semistable.
\end{proposition}
\begin{proof}
Assume $(n,d) \neq (5,2)$. Consider the set
\begin{equation*}
I':=\big\{{X_0}^d,\,{X_0}^{d-1}X_1,\ldots,X_0{X_1}^{d-1},
    \,{X_1}^d,\,{X_2}^d\big\}.
\end{equation*}
and the sequence
\begin{align*}
\big(&X_0{X_2}^{d-1},\,{X_0}^2{X_2}^{d-2},\ldots,{X_0}^{d-2}{X_2}^2,\\
    &{X_1}^{d-1}X_2,\,{X_0}^{d-1}X_2,\\
    &X_1{X_2}^{d-1},\,{X_1}^2{X_2}^{d-2},\ldots,{X_1}^{d-2}{X_2}^2\big).
\end{align*}
If $1\leq i\leq 2d-2$, let $I''$ be the set of the first~$i$ monomials in this sequence and let ${I_{d,n}=I'\cup I''}$. The number of monomials in~$I_{d,n}$ is ${n=d+2+i}$.

For ${0<d_J<d}$, since we are again following the strategy mentioned in Proposition~\ref{interval2}, it is enough to count, in each case, the number of multiples of ${X_0}^{d_J}$ which are in $I_{d,n}$.

If $i\leq d-2$, the set of multiples of ${X_0}^{d_J}$ in~$I_{d,n}$ is
\[\big\{{X_0}^d,\ldots,{X_0}^{d_J}{X_1}^{d-d_J},\,
    {X_0}^{d_J}{X_2}^{d-d_J},\ldots,{X_0}^e{X_2}^{d-e}\big\},\]
where $e:=\max\{i,d_J-1\}$ and the list ${{X_0}^{d_J}{X_2}^{d-d_J},\ldots,{X_0}^e{X_2}^{d-e}}$ is understood to be empty if ${e=d_J-1}$. The number of monomials in this set is ${k=d-2d_J+e+2}$, and we get
\begin{align*}
(d-d_J)n+d_J-dk&=i(d-d_J)+dd_J-d_J-de>0.
\end{align*}
If $i=d-1$, the set of multiples of~${X_0}^{d_J}$ is
\[\big\{{X_0}^d,\ldots,{X_0}^{d_J}{X_1}^{d-d_J},\,
    {X_0}^{d_J}{X_2}^{d-d_J},\ldots,{X_0}^{d-2}{X_2}^2\big\}.\]
The list ${X_0}^{d_J}{X_2}^{d-d_J},\ldots,{X_0}^{d-2}{X_2}^2$ is again understood to be empty if ${d_J=d-1}$. The number of monomials in this set is ${k=2d-2d_J}$, and we get
\begin{align*}
(d-d_J)n+d_J-dk&=d>0.
\end{align*}
If $i\geq d$, the set of multiples of~${X_0}^{d_J}$ is
\[\big\{{X_0}^d,\ldots,{X_0}^{d_J}{X_1}^{d-d_J},\,
    {X_0}^{d_J}{X_2}^{d-d_J},\ldots,{X_0}^{d-1}X_2\big\}.\]
The number of monomials in this set is ${k=2d-2d_J+1}$, and we get
\begin{align*}
(d-d_J)n+d_J-dk&\geq d-d_J>0.
\end{align*}

In all cases, inequality~(\ref{ineq4}) is strictly satisfied, and the corresponding  syzygy bundle is stable.

For $(n,d)=(5,2)$, it is enough to take
$I=\{{X_0}^2,\,{X_1}^2,\,{X_2}^2,\,X_0X_1,\,X_0X_2\}$.
\end{proof}

\begin{proposition}\label{interval4} For any integers $n$ and $d$ such that  $3d<n\le\tbinom{d+2}{2}$, there is a set $I_{d,n}$ of $n$  $\frak{m}$-primary monomials in $K[X_0,X_1,X_2]$ of degree $d$ such that the corresponding syzygy bundle $E_{d,n}$ is stable.
\end{proposition}
\begin{proof}
We divide the proof in three cases. Let $j\geq1$ be such that $3j<d$ and suppose that
\[\tbinom{d+2}{2}-\tbinom{d+2-3j}{2}<n\leq
    \tbinom{d+2}{2}-\tbinom{d+2-3(j+1)}{2}.\]
Note that as $j$ varies, we get all values of~$n$ mentioned, except $\tbinom{d+2}{2}$ when~$d$ is a multiple of~$3$. However, for this highest possible value of~$n$, the result follows from Proposition \ref{homogeni}.

\bigskip

\noindent\textbf{Case~1.}
Suppose that
\[{n=\tbinom{d+2}{2}-\tbinom{d+2-3j}{2}+i=
    3dj-\tfrac{9j(j-1)}{2}+i},\]
with ${1\leq i\leq d-3j+1}$ and consider the set
\[I':=\big\{{X_0}^{i_0}{X_1}^{i_1}{X_2}^{i_2}:i_0+i_1+i_2=d \text{ and }
    (i_0<j\lor i_1<j\lor i_2<j)\big\}.\]
Consider the sequence
\begin{align*}
\big(&{X_0}^{d-2j}{X_1}^j{X_2}^j,\,
    {X_0}^{d-2j-1}{X_1}^{j+1}{X_2}^j,
    \ldots,{X_0}^j{X_1}^{d-2j}{X_2}^j\big).
\end{align*}
Let $I''$ be the set of the first~$i$ monomials in this sequence and let ${I_{d,n}=I'\cup I''}$. Then~$I_{d,n}$ has~$n$ monomials and we  verify that it strictly satisfies inequality~(\ref{ineq4}).
\[
\xymatrix@C=\meiabase@R=\altura@!0{
&&&&&&&&&&&&{\bullet}\\
&&&&&&&&&&&{\bullet}&&{\bullet}\\
&&&&&&&&&&{\bullet}&&{\bullet}&&{\bullet}\\
&&&&&&&&&{\bullet}&&{\bullet}&&{\bullet}&&{\bullet}\\
&&&&&&&&{\bullet}&&{\bullet}&&{\circ}&&{\bullet}&&{\bullet}\\
&&&&&&&{\bullet}&&{\bullet}&&{\circ}&&{\circ}&&{\bullet}&&{\bullet}\\
&&&&&&{\bullet}&&{\bullet}&&{\circ}&&{\circ}&&{\circ}&&{\bullet}&&
    {\bullet}\\
&&&&&{\bullet}&&{\bullet}&&{\circ}&&{\circ}&&{\circ}&&{\circ}&&
    {\bullet}&&{\bullet}\\
&&&&{\bullet}&&{\bullet}&&{\circ}&&{\circ}&&{\circ}&&{\circ}&&
    {\circ}&&{\bullet}&&{\bullet}\\
&&&{\bullet}&&{\bullet}&&{\circ}&&{\circ}&&{\circ}&&{\circ}&&
    {\circ}&&{\circ}&&{\bullet}&&{\bullet}\\
&&{\bullet}&&{\bullet}&&{\bullet}&&{\bullet}&&{\bullet}&&{\circ}&&
    {\circ}&&{\circ}&&{\circ}&&{\bullet}&&{\bullet}\\
&{\bullet}&&{\bullet}&&{\bullet}&&{\bullet}&&{\bullet}&&{\bullet}&&
    {\bullet}&&{\bullet}&&{\bullet}&&{\bullet}&&{\bullet}&&
    {\bullet}\\
{\bullet}&&{\bullet}&&{\bullet}&&{\bullet}&&{\bullet}&&{\bullet}&&
    {\bullet}&&{\bullet}&&{\bullet}&&{\bullet}&&{\bullet}&&
    {\bullet}&&{\bullet}\\
\\
&&&&&&&&&&&&I_{12,66}
}
\]

For ${0<d_J<d}$, all we have to do is to count, in each case the number of multiples of~${X_0}^{d_J}$ which are present in~$I_{d,n}$, since we are again applying the strategy mentioned in Proposition~\ref{interval2}.

For $d-2j\leq d_J<d$, all monomials of degree~$d$ of type ${X_0}^{i_0}{X_1}^{i_1}{X_2}^{i_2}$, with $i_0\geq d_J$, are in~$I_{d,n}$. Therefore the number of multiples of~${X_0}^{d_J}$ in~$I_{d,n}$ is
\[k=\tbinom{d-d_J+2}{2}\]
and we get
\begin{align*}
(d-d_J)n+d_J-dk&=
    (d-d_J)\left(3dj-\tfrac{9j(j-1)}{2}+i\right)
        +d_J-d\tbinom{d-d_J+2}{2}.
\end{align*}
This expression can be rewritten in the two following ways:
\begin{align*}
&\tfrac{1}{2}d(d-d_J)(d_J+j-d)
        +\tfrac{5}{2}(d-3j)(d-d_J)(j-1)
        + 3(d-d_J)j(j-1) +\\
    &\qquad +d(d-d_J) + (i-1)(d-d_J)\\
\intertext{and}
&\tfrac{1}{2}d(d-d_J-j)(d_J+2j-d)
        + \tfrac{3}{2}(d-d_J-j)^2(j-1)+\\
    &\qquad+ \tfrac{3}{2}(d-d_J-j)d_J(j-1)
        + \tfrac{1}{2}(d_J+2j-d)j^2+\\
    &\qquad+ \tfrac{5}{2}(d_J-j)j(j-1)
        + \tfrac{3}{2}(d-d_J-j)j + d_Jj
        + \tfrac{1}{2}j^2 + (i-1)(d-d_J).
\end{align*}
>From the first one, we can see that the expression above is positive for ${d-j\le d_J<d}$, and the second shows us positivity for ${d-2j\le d_J<d-j}$ (since ${3j<d}$, we get in this case ${j<d_J}$).

For $j\leq d_J<d-2j$, the monomials in~$I_{d,n}$ that are multiples of~${X_0}^{d_J}$ are the ones in the set
\begin{multline*}
J:=\big\{{X_0}^{i_0}{X_1}^{i_1}{X_2}^{i_2}\in I':
    i_0\geq d_J\big\}\cup\\
    \cup\big\{{X_0}^{d-2j}{X_1}^j{X_2}^j,\,
    {X_0}^{d-2j-1}{X_1}^{j+1}{X_2}^j,
    \ldots,{X_0}^{d-2j-e}{X_1}^{j+e}{X_2}^j\big\},
\end{multline*}
where $e:=\min(i-1,d-2j-d_J)$. Therefore their number is \[{k=\tbinom{d-d_J+2}{2}-\tbinom{d-2j-d_J+2}{2}+e}.\]
If $i-1\leq d-2j-d_J$, we get
\begin{align*}
(d-d_J)n+d_J-dk&=
    (d-d_J)\left(3dj-\tfrac{9j(j-1)}{2}+i\right)+d_J-\\
    &\quad-d\left(\tbinom{d-d_J+2}{2}
        -\tbinom{d-2j-d_J+2}{2}+i-1\right)=\\
    &=(d-j-d_J)(d-2j-d_J)j + (d-2j-d_J)d_J(j-1) +\\
    &\qquad + \tfrac{1}{2}dj^2 + \tfrac{7}{2}(d_J-j)j(j-1)
        + \tfrac{3}{2}j^2(j-1) +\\
    &\qquad + \tfrac{3}{2}(d-2j-d_J)j + \tfrac{1}{2}d_Jj
        + j^2 + d +\\
    &\qquad + (d-2j-d_J+1-i)d_J>0
\end{align*}
since $d-2j-d_J>0$ and $j\geq1$. If $i-1>d-2j-d_J$, we get
\begin{align*}
(d-d_J)n+d_J-dk&=
    (d-d_J)\left(3dj-\tfrac{9j(j-1)}{2}+i\right)+d_J-\\
    &\quad-d\left(\tbinom{d-d_J+2}{2}
        -\tbinom{d-2j-d_J+2}{2}+d-2j-d_J\right)=\\
    &=(d-d_J-2j)(d-j)(j-1) + \tfrac{1}{2}dj^2 +
        \tfrac{7}{2}(d_J-j)j(j-1) +\\
    &\quad +\tfrac{3}{2}j^3 + (d-d_J-2j)^2 +
        \tfrac{5}{2}(d-d_J-2j)j
        + \tfrac{1}{2}(d_J-j)j + d +\\
    &\quad + (i-1-d+2j+d_J)(d-d_J)>0
\end{align*}
since $d-2j-d_J>0$ and $j\geq1$.

For $0<d_J<j$, the number of monomials in~$I_{d,n}$ that are multiples of~${X_0}^{d_J}$ is \[{k=\tbinom{d-d_J+2}{2}-\tbinom{d-3j+2}{2}+i}\]
and we get
\begin{align*}
(d-d_J)n+d_J-dk&=
    (d-d_J)\left(3dj-\tfrac{9j(j-1)}{2}+i\right)+d_J-\\
    &\quad-d\left(\tbinom{d-d_J+2}{2}
        -\tbinom{d-3j+2}{2}+i\right)=\\
    &=(d-3j)^2d_J + 3(d-3j)(j-d_J)d_J
        + \tfrac{5}{2}(d-3j)d_J^2+\\
    &\quad + \tfrac{9}{2}j(j-d_J)d_J + 3jd_J^2
        + \tfrac{1}{2}(d-3j)d_J+ (d-3j+1-i)d_J>0.
\end{align*}
In all cases, inequality~(\ref{ineq4}) is strictly satisfied, and hence the corresponding syzygy bundle is stable.

\bigskip

\noindent\textbf{Case~2.}
Now suppose that
\[n=\tbinom{d+2}{2}-\tbinom{d+1-3j}{2}+i=
    3dj+d+1-\tfrac{3j(3j-1)}{2}+i,\]
with ${1\leq i\leq d-3j}$ and consider the set
\[I':=\big\{{X_0}^{i_0}{X_1}^{i_1}{X_2}^{i_2}:i_0+i_1+i_2=d \text{ and }
    (i_0<j\lor i_1<j\lor i_2\leq j)\big\}.\]
Consider the sequence
\begin{align*}
\big(&{X_0}^j{X_1}^j{X_2}^{d-2j},\,
    {X_0}^{j+1}{X_1}^j{X_2}^{d-2j-1},
    \ldots,{X_0}^{d-2j-1}{X_1}^j{X_2}^{j+1}\big).
\end{align*}
Let $I''$ be the set of the first~$i$ monomials in this sequence and let ${I_{d,n}=I'\cup I''}$. Then~$I_{d,n}$ has~$n$ monomials and we  verify that it strictly satisfies inequality~(\ref{ineq4}).
\[
\xymatrix@C=\meiabase@R=\altura@!0{
&&&&&&&&&&&&{\bullet}\\
&&&&&&&&&&&{\bullet}&&{\bullet}\\
&&&&&&&&&&{\bullet}&&{\bullet}&&{\bullet}\\
&&&&&&&&&{\bullet}&&{\bullet}&&{\bullet}&&{\bullet}\\
&&&&&&&&{\bullet}&&{\bullet}&&{\bullet}&&{\bullet}&&{\bullet}\\
&&&&&&&{\bullet}&&{\bullet}&&{\bullet}&&{\circ}&&{\bullet}&&{\bullet}\\
&&&&&&{\bullet}&&{\bullet}&&{\bullet}&&{\circ}&&{\circ}&&{\bullet}&&
    {\bullet}\\
&&&&&{\bullet}&&{\bullet}&&{\circ}&&{\circ}&&{\circ}&&{\circ}&&
    {\bullet}&&{\bullet}\\
&&&&{\bullet}&&{\bullet}&&{\circ}&&{\circ}&&{\circ}&&{\circ}&&
    {\circ}&&{\bullet}&&{\bullet}\\
&&&{\bullet}&&{\bullet}&&{\circ}&&{\circ}&&{\circ}&&{\circ}&&
    {\circ}&&{\circ}&&{\bullet}&&{\bullet}\\
&&{\bullet}&&{\bullet}&&{\bullet}&&{\bullet}&&{\bullet}&&{\bullet}&&
    {\bullet}&&{\bullet}&&{\bullet}&&{\bullet}&&{\bullet}\\
&{\bullet}&&{\bullet}&&{\bullet}&&{\bullet}&&{\bullet}&&{\bullet}&&
    {\bullet}&&{\bullet}&&{\bullet}&&{\bullet}&&{\bullet}&&
    {\bullet}\\
{\bullet}&&{\bullet}&&{\bullet}&&{\bullet}&&{\bullet}&&{\bullet}&&
    {\bullet}&&{\bullet}&&{\bullet}&&{\bullet}&&{\bullet}&&
    {\bullet}&&{\bullet}\\
\\
&&&&&&&&&&&&I_{12,73}
}
\]

As in the previous step, for ${0<d_J<d-1}$, no monomial of degree~$d_J$ divides a greater number of monomials in~$I_{d,n}$ than~${X_0}^{d_J}$. Therefore all we have to do is count, in each case the number of multiples of~${X_0}^{d_J}$ which are present in~$I_{d,n}$.

For $d-2j\leq d_J<d$, all monomials of degree~$d$ of type ${X_0}^{i_0}{X_1}^{i_1}{X_2}^{i_2}$, with $i_0\geq d_J$,
are in~$I_{d,n}$. Therefore the number of multiples of~${X_0}^{d_J}$ in~$I_{d,n}$ is
\[k=\tbinom{d-d_J+2}{2},\]
as it was in step~1, and we can claim that since all values are
the same except for $n$, which is bigger, inequality~(\ref{ineq3})
is strictly satisfied, due to the fact that sequence
$(a_{d,j})_{j\geq 2}$ is monotonically increasing (see Remark
\ref{rem}(b)).

For $j\leq d_J<d-2j$, the monomials in~$I_{d,n}$ that are multiples of~${X_0}^{d_J}$ are the ones in the set
\begin{multline*}
J:=\big\{{X_0}^{i_0}{X_1}^{i_1}{X_2}^{i_2}\in I':
    i_0\geq d_J\big\}\cup\\
    \cup\big\{{X_0}^{d_J}{X_1}^j{X_2}^{d-j-d_J},\,
    {X_0}^{d_J+1}{X_1}^j{X_2}^{d-j-d_J-1},
    \ldots,{X_0}^{j+i-1}{X_1}^j{X_2}^{d-2j-i+1}\big\},
\end{multline*}
where this last set is understood to be empty if ${j+i-1<d_J}$.
Therefore their number is
\[k=\tbinom{d-d_J+2}{2}-\tbinom{d-2j-d_J+1}{2}+\max(0,j+i-d_J).\]
If $j+i\le d_J$, we get (keeping in mind that $i\geq1$)
\begin{align*}
(d-d_J)n+d_J-dk&=
    (d-d_J)\left(3dj+d+1-\tfrac{3j(3j-1)}{2}+i\right)+d_J-\\
    &\quad-d\left(\tbinom{d-d_J+2}{2}-\tbinom{d-2j-d_J+1}{2}
        \right)=\\
    &=(d-j)(d-2j-d_J)j + \tfrac{1}{2}dj^2
        + \tfrac{7}{2}(d_J-j)j(j-1) + \tfrac{3}{2}j^3+\\
    &\quad + \tfrac{1}{2}(d-2j-d_J)j + \tfrac{5}{2}(d_J-j)j
        + i(d-d_J)>0.
\end{align*}
If $j+i>d_J$, we get
\begin{align*}
(d-d_J)n+d_J-dk&=
    (d-d_J)\left(3dj+d+1-\tfrac{3j(3j-1)}{2}+i\right)+d_J-\\
    &\quad-d\left(\tbinom{d-d_J+2}{2}-\tbinom{d-2j-d_J+1}{2}
        +j+i-d_J\right)=\\
    &=(d-j)(d-2j-d_J)j + \tfrac{1}{2}(d-2j-d_J)j(j-1)+\\
    &\quad + 4(d_J-j)j^2 + 3j^3+ (d_J-j)j+ (d-3j-i)d_J>0.
\end{align*}

For $0<d_J<j$, the number of monomials in~$I_{d,n}$ that are multiples of~${X_0}^{d_J}$ is
\[k=\tbinom{d-d_J+2}{2}-\tbinom{d-3j+1}{2}+i\]
and we get
\begin{align*}
(d-d_J)n+d_J-dk&=
    (d-d_J)\left(3dj+d+1-\tfrac{3j(3j-1)}{2}+i\right)+d_J-\\
    &\quad-d\left(\tbinom{d-d_J+2}{2}-\tbinom{d-3j+1}{2}+i
        \right)=\\
    &=(d-j)(d-3j)d_J + d(j-d_J)d_J
        +\tfrac{1}{2}dd_J(d_J-1)+\\
    &\quad + \tfrac{3}{2}j^2d_J
        + \tfrac{3}{2}jd_J+ (d-3j-i)d_J>0.
\end{align*}
Again in all cases, inequality~(\ref{ineq4}) are strictly satisfied and the associated  syzygy bundle is stable.

\bigskip

\noindent\textbf{Case~3.}
If $d=3j+1$, case~2 has exhausted all possible monic monomials of degree~$d$, and this proof is ended.

If ${d>3j+1}$, then suppose that
\[n=\tbinom{d+2}{2}-\tbinom{d-3j}{2}+i=
    3dj+2d+1-\tfrac{3j(3j+1)}{2}+i,\]
with ${1\leq i\leq d-3j-1}$, and consider the set
\[I':=\big\{{X_0}^{i_0}{X_1}^{i_1}{X_2}^{i_2}:i_0+i_1+i_2=d
    \text{ and }
    (i_0<j\lor i_1\leq j\lor i_2\leq j)\big\}.\]
Consider the ordered multiple
\begin{align*}
\big(&{X_0}^j{X_1}^{j+1}{X_2}^{d-2j-1},\,
    {X_0}^j{X_1}^{j+2}{X_2}^{d-2j-2},
    \ldots,{X_0}^j{X_1}^{d-2j-1}{X_2}^{j+1}\big).
\end{align*}
Let $I''$ be the set of the first~$i$ monomials in this ordered multiple and let ${I_{d,n}=I'\cup I''}$. Then~$I_{d,n}$ has~$n$ monomials and we shall verify that it strictly satisfies inequality~(\ref{ineq4}).
\[
\xymatrix@C=\meiabase@R=\altura@!0{
&&&&&&&&&&&&{\bullet}\\
&&&&&&&&&&&{\bullet}&&{\bullet}\\
&&&&&&&&&&{\bullet}&&{\bullet}&&{\bullet}\\
&&&&&&&&&{\bullet}&&{\bullet}&&{\bullet}&&{\bullet}\\
&&&&&&&&{\bullet}&&{\bullet}&&{\bullet}&&{\bullet}&&{\bullet}\\
&&&&&&&{\bullet}&&{\bullet}&&{\bullet}&&{\bullet}&&{\bullet}&&{\bullet}\\
&&&&&&{\bullet}&&{\bullet}&&{\bullet}&&{\circ}&&{\bullet}&&{\bullet}&&
    {\bullet}\\
&&&&&{\bullet}&&{\bullet}&&{\bullet}&&{\circ}&&{\circ}&&{\circ}&&
    {\bullet}&&{\bullet}\\
&&&&{\bullet}&&{\bullet}&&{\bullet}&&{\circ}&&{\circ}&&{\circ}&&
    {\circ}&&{\bullet}&&{\bullet}\\
&&&{\bullet}&&{\bullet}&&{\bullet}&&{\circ}&&{\circ}&&{\circ}&&
    {\circ}&&{\circ}&&{\bullet}&&{\bullet}\\
&&{\bullet}&&{\bullet}&&{\bullet}&&{\bullet}&&{\bullet}&&{\bullet}&&
    {\bullet}&&{\bullet}&&{\bullet}&&{\bullet}&&{\bullet}\\
&{\bullet}&&{\bullet}&&{\bullet}&&{\bullet}&&{\bullet}&&{\bullet}&&
    {\bullet}&&{\bullet}&&{\bullet}&&{\bullet}&&{\bullet}&&
    {\bullet}\\
{\bullet}&&{\bullet}&&{\bullet}&&{\bullet}&&{\bullet}&&{\bullet}&&
    {\bullet}&&{\bullet}&&{\bullet}&&{\bullet}&&{\bullet}&&
    {\bullet}&&{\bullet}\\
\\
&&&&&&&&&&&&I_{12,78}
}
\]

As in the previous cases, all we have to do is count, in each case the number of multiples of~${X_0}^{d_J}$ which are present in~$I_{d,n}$.

For $d-2j\leq d_J<d$, all monomials of degree~$d$ of type ${X_0}^{i_0}{X_1}^{i_1}{X_2}^{i_2}$, with $i_0\geq d_J$, are in~$I_{d,n}$. Therefore the number of multiples of~${X_0}^{d_J}$ in~$I_{d,n}$ is
\[k=\tbinom{d-d_J+2}{2},\]
as it was in cases~1 and~2, and we can claim that since all values
are the same except for~$n$, which is bigger,
inequality~(\ref{ineq4}) are strictly satisfied, due to the fact
that sequence $(a_{d,j})_{j\geq2}$ is monotonically increasing
(see Remark \ref{rem}(b)).

For $j<d_J<d-2j$, an analogous argument based on calculations for case~2 allows us to claim that inequality~(\ref{ineq4}) is strictly satisfied.

For $0<d_J\leq j$, the number of monomials in~$I_{d,n}$ that are multiples of~${X_0}^{d_J}$ is
\[k=\tbinom{d-d_J+2}{2}-\tbinom{d-3j}{2}+i\]
and we get (keeping in mind that $i\leq d-3j-1$)
\begin{align*}
(d-d_J)n+d_J-dk&=
    (d-d_J)\left(3dj+2d+1-\tfrac{3j(3j+1)}{2}+i\right)+d_J-\\
    &\quad-d\left(\tbinom{d-d_J+2}{2}-\tbinom{d-3j}{2}+i
        \right)=\\
    &=(d-2j)(d-3j)d_J + 2(d-j)(j-d_J)d_J
        + \tfrac{3}{2}(d-j)d_J(d_J-1)+\\
    &\quad + \tfrac{1}{2}j(j-d_J)d_J
        + 3jd_J + d_J + (d-3j-1-i)d_J>0.
\end{align*}
Again in all cases, inequality~(\ref{ineq4}) is strictly satisfied, which makes the syzygy bundle stable, and concludes the proof.
\end{proof}

Putting all together we have got

\begin{theorem} \label{mainthm}
For any integers $d,n \ge 1$ with $(n,d) \neq (5,2)$ and $3\le n
\le\tbinom{d+2}{2}$, there is a family of $n$ $\frak{m}$-primary
monomials in $K[X_0,X_1,X_2]$ of degree $d$ such that the
corresponding syzygy bundle is stable. For $(n,d) =(5,2)$, there
are $5$ $\frak{m}$-primary monomials in $K[X_0,X_1,X_2]$ of degree
$2$ such that the corresponding syzygy bundle is semistable.
\end{theorem}

As an immediate consequence  of Theorem \ref{mainthm} we obtain

\begin{corollary} Let $E_{d,n}$ be the syzygy bundle on $\PP^2$ associated to $n$ general
$\frak{m}$-primary forms of the same degree $d$. Suppose that
$3\le n\le\tbinom{d+2}{2}$. Then $E_{d,n}$ is stable when $(n,d)
\neq (5,2)$ and $E_{2,5}$ is semistable.
\end{corollary}
\begin{proof}
It follows from Theorem \ref{mainthm}, taking into account that stability is an open property.
\end{proof}

\section{Moduli spaces of syzygy bundles}

In this section we  study the moduli space of syzygy bundles
on $\PP^N$. We denote by $M=M(r;c_1, \ldots, c_{s})$ the moduli space of rank $r$, stable vector bundles $E$ on  $\PP^N$ with fixed Chern classes $c_i(E)=c_i$, for $i=1,\ldots,s=min(r,N)$. The existence of the moduli space $M(r;c_1, \ldots, c_{s})$ was established by Maruyama in 1977 (see \cite{MA} and \cite{MA1}) and once the existence of the moduli space is established, the question arises as what can be said about its local and global structure. More precisely, what does the moduli space look like as an algebraic variety? Is it, for example, connected, irreducible, rational or smooth? What does it look like as a topological space? What is its geometry?  Until now, there is no general answer to these questions. The goal of this section is to determine the unobstructedness of stable syzygy bundles $E_{d,n}$ on $\PP^N$ and to compute the dimension of the irreducible component of the corresponding moduli space.

Let us start by analyzing whether a syzygy bundle on $\PP^N$ is stable and to state our contribution to study (semi)stability
properties of syzygy bundles on $\PP^N$. This will improve all previous known results, which we quickly recall now.

\vskip 2mm Let $C\subset \PP^N$ be a smooth, projective, elliptic curve embedded by a complete system of degree $N+1$. Using the fact that the restriction of a general syzygy bundle $E_{d,n}$ on $\PP^N$ to $C$ is (semi)stable, Hein proved:

\begin{proposition}
Let $f_1,\ldots,f_n\in K[X_0,X_1,\ldots,X_N]$, $N\ge 2$, denote generic homogeneous forms of degree $d$. Suppose that
$N+1\le n\le d(N+1)$. Then the syzygy bundle $E_{d,n}$ on $\PP^N$ is semistable.
\end{proposition}
\begin{proof}
See \cite{B}, Theorem 8.6 and Theorem A.1.
\end{proof}

\vskip 2mm
As another application of Theorem \ref{mainthm}, we can improve the above proposition and we get

\begin{theorem} \label{mainthm2}
Let $f_1,\ldots ,f_n\in K[X_0,X_1,\ldots ,X_N]$, $N\ge 2$, denote generic homogeneous polynomials of degree $d$. Suppose $N+1\le n\le\tbinom{d+2}{2}+N-2$. Then the syzygy bundle $E_{d,n}$ on $\PP^N$ is stable  when ${(N,n,d) \neq (2,5,2)}$, and $E_{2,5}$ is semistable on $\PP^2$.
\end{theorem}
\begin{proof}
Since stability is an open property in a flat family of torsion free sheaves, it is enough to prove the stability property for a single choice of homogeneous forms $f_1,\ldots ,f_n$ of degree $d$.

If $(n,d) \neq (5,2)$ we proceed by induction on $N$. By induction on the number of variables, the case $N=2$ being done, we can suppose that for given $N$, $n$ and $d$ in the above conditions, there is a set of primary monomials $I$ such that for any $J\subset I$, with $|J|\ge 2$, inequality (\ref{ineq3}) is valid.

Let $n$ be an integer such that $N+2\le n\le \tbinom{d+2}{2}+(N+1)-2$. Since
\[N+1\le n-1\le \tbinom{d+2}{2}+N-2,\]
there is a set $I_0$ of $n-1$ primary monomials in $K[X_0,\cdots ,X_N]$ satisfying inequality (\ref{ineq3}). Let
\[I:=I_0\cup \big\{X_{N+1}^d\big\}.\]
Then the ideal generated by $I$ is primary. Let $J\subset I$ be a subset with at least two monomilas. If $J\subset I_0$, then by induction hypothesis, inequality (\ref{ineq3}) holds. If not, then ${X_{N+1}}^d \in J$, and since $J$ has at least another monomial, where the variable $X_{N+1}$ does not occur, the greatest common divisor is 1, and $d_J=0$; so inequality (\ref{ineq3}) holds because the sequence $\left(a_{d,j}\right)_{j\ge 2}$ is monotonically increasing.

Assume $(n,d)=(5,2)$. Note that in that case $2 \le N \le 4$. If $N=2$, $E_{5,2}$ is a semistable bundle on $\PP^2$ by Theorem~\ref{mainthm}. If $N=4$, $E_{5,2}$ is a stable bundle on $\PP^4$ by Proposition~\ref{BSpindler}. Finally if $N=3$, we consider the set $I:=\big\{{X_0}^2,\,{X_1}^2,\,{X_2}^2,\,{X_3}^2,\,X_0X_1\big\}$. The associated syzygy bundle is stable and hence, by the openness of the stability, $E_{5,2}$ is stable on $\PP^3$.
\end{proof}

\begin{remark}In general, the bound $N+1\le n\le\tbinom{d+2}{2}+N-2$  generalizes the bound $N+1 \le n \le d(N+1)$ given by Hein in \cite{B1}, Theorem A.1.
\end{remark}

We are now ready to state the unobstructedness of stable syzygy bundles on $\PP^N$.

\begin{theorem}\label{moduli}
Assume $N+1\le n\le\tbinom{d+2}{2}+N-2$, $N\ne 3$ and $(N,n,d)
\neq (2,5,2)$. Then the syzygy bundle $E_{d,n}$ is unobstructed
and it belongs to a generically smooth irreducible component of
dimension $n\tbinom{d+N}{N}-n^2$, if $N \geq 4$, and
$n\tbinom{d+2}{2}+n\tbinom{d-1}{2}-n^2$, if $N=2$.
\end{theorem}
\begin{proof}
Let us denote by $c_i=c_i(E_{d,n})$, $i=1,\ldots,\min(n-1, N)$ the $i$th Chern class of $E_{d,n}$ and let $M=M(n-1;c_1,\ldots,c_{\min(n-1,N)})$ be the moduli space of rank $n-1$, stable vector bundles  on  $\PP^N$ with Chern classes $c_i$. From deformation theory, we know that the Zariski tangent space of $M$ at $[E_{d,n}]$ is canonically given by
\[T_{[E_{d,n}]}M\cong Ext^{1}(E_{d,n},E_{d,n})
    \cong H^1(E_{d,n} \otimes E_{d,n}^{\vee});\]
and the obstruction space of the local ring $\cO_{M,[E_{d,n}]}$ is a subspace of $Ext^{2}(E_{d,n},E_{d,n})$. Thus, if
\[Ext^2(E_{d,n},E_{d,n}) \cong H^2(E_{d,n} \otimes
    E_{d,n}^{\vee})=0,\]
then the moduli space $M$ is smooth at $E_{d,n}$ and in this last case
\[\dim_{K}Ext^{1}(E_{d,n},E_{d,n})=
    \dim_{[E_{d,n}]}M(n-1;c_1,\ldots,c_{\min(n-1,N)})\]
(see \cite{MA} and \cite{MA1}).

To compute $Ext^{i}(E_{d,n},E_{d,n})$, we consider the exact sequence
\begin{equation} \label{suc1}
\xymatrix{
0\ar[r] &E_{d,n}\ar[r] &\cO_{\PP^N}(-d)^{n}\ar[r]
    &\cO_{\PP^N}\ar[r] &0
}
\end{equation}
and its dual
\begin{equation}\label{suc2}
\xymatrix{
0\ar[r] &\cO_{\PP^N}\ar[r] &\cO_{\PP^N}(d)^{n}\ar[r]
    &E_{d,n}^{\vee}\ar[r] &0.
}
\end{equation}

First of all, note that by the cohomological exact sequence
associated to the exact sequence (\ref{suc1}), we get
\begin{equation} \label{coho1}
\begin{array}{l} h^0(E_{d,n})=0; \\
h^1(E_{d,n})=1; \\
h^2(E_{d,n})=
\left\{
    \begin{array}{ll}
        0, & \mbox{if} \quad N \geq 4 \\
        n\tbinom{d-1}{2}, & \mbox{if} \quad N=2;
    \end{array}
\right.\\
h^3(E_{d,n})=0.
\end{array}
\end{equation}

Denote by $F =E_{d,n} \otimes E_{d,n}^{\vee}$. Consider the
cohomological exact sequence
\begin{equation}\label{ext2}
\begin{split}
\xymatrix@R=0pt{
0\ar[r] &H^0(E_{d,n})\ar[r] &H^0(E_{d,n}(d))^n
    \ar[r] &H^0(F)\ar[r] &\\
\ar[r] &H^1(E_{d,n})\ar[r] &H^1(E_{d,n}(d))^n\ar[r]
    &H^1(F)\ar[r] &\\
\ar[r] &H^2(E_{d,n})\ar[r] &H^2(E_{d,n}(d))^n\ar[r]
    &H^2(F)\ar[r] &\\
\ar[r] &H^3(E_{d,n})\ar[r] &H^3(E_{d,n}(d))^n\ar[r]
    &H^3(F)\ar[r] &\cdots
}
\end{split}
\end{equation}
associated to the exact sequence
\[
\xymatrix{
0\ar[r] &E_{d,n}\ar[r] &E_{d,n}(d)^{n}\ar[r] &F\ar[r] &0.
}
\]
Since $E_{d,n}$ is stable, it is simple, i.e. $H^0(F)=K$. Thus, from the exact sequence (\ref{ext2}), and the fact that by (\ref{coho1}), $H^0(E_{d,n})=0$, we get $H^0(E_{d,n}(d))=0$.

Twisting by $\cO_{\PP^N}(d)$ the exact sequence (\ref{suc1}), and taking cohomology, we deduce
\begin{equation} \label{coho2}
\begin{split}
    &h^2(E_{d,n}(d))=0, \\
    &h^3(E_{d,n}(d))=0, \\
    &h^1(E_{d,n}(d))=\tbinom{N+d}{d}-n.
\end{split}
\end{equation}

In particular, from (\ref{ext2}) we get $H^2(F) \cong
Ext^2(E_{d,n},E_{d,n})=0$ and the exact sequence
\[
\xymatrix{
0\ar[r] &K\ar[r] &K\ar[r] &H^1(E_{d,n}(d))^n\ar[r]
    &H^1(F)\ar[r] &H^2(E_{d,n})\ar[r] &0.
}
\]
Therefore
\[ h^1(F) = ext^1(E_{d,n},E_{d,n})= \left\{
    \begin{array}{ll}
        n\tbinom{N+d}{d}-n^2, & \mbox{if} \quad N \ge4\\
        n\tbinom{d+2}{2}+n\tbinom{d-1}{2}-n^2, & \mbox{if}
            \quad N=2,
   \end{array}
\right.
\]
which finishes the proof.
\end{proof}

\end{document}